\documentclass[a4paper,oneside,11pt]{article}

\usepackage[latin1]{inputenc}
\usepackage[T1]{fontenc}
\usepackage[english,french]{babel}
\usepackage[a4paper,vmargin={2.5cm,3.5cm},hmargin=2cm]{geometry}
\usepackage[colorlinks=false]{hyperref}
\usepackage[expansion=false]{microtype}
\usepackage{graphicx}
\usepackage{xcolor}

\usepackage{amsmath,amsfonts,amssymb,amsthm,mathrsfs}

\theoremstyle{plain}
\newtheorem{theoreme}{Th\'eor\`eme}
\newtheorem{corollaire}{Corollaire}
\newtheorem{proposition}{Proposition}
\newtheorem{lemme}{Lemme}
\theoremstyle{definition}
\newtheorem{definition}{D\'efinition}
\newtheorem{etape}{\'Etape}
\theoremstyle{remark}
\newtheorem*{remarque}{Remarque}

\hypersetup{
    pdfauthor   = {Laurent Menard},%
    pdftitle    = {La quadrangulation infinie uniforme}}

\author{Laurent {\sc M\'enard}\footnote{\texttt{laurent.menard@normalesup.org}}\\
        D\'epartement de Math\'ematiques\\
        Universit\'e Paris-Sud\\
        91405 Orsay Cedex, France}

\title{\bf LES DEUX QUADRANGULATIONS INFINIES UNIFORMES ONT M\^EME LOI}
\date{{\small May 2008}}

\newcommand*\bij{\Phi}
\newcommand*\arbres{\overline{\mathbb{T}}}
\newcommand*\arbre[1][\infty]{\mathbb{T}_{#1}}
\newcommand*\spine{\mathcal{S}}
\newcommand*\spinelab{\mathcal{C}}
\newcommand*\carte[1][\infty]{\mathbf{Q}_{#1}}
\newcommand*\cartes{\overline{\mathbb{Q}}}
\newcommand*\kcartes{\overline{\mathbf{Q}}}
\newcommand*\boule[2][R]{B_{#2,#1}}
\newcommand*\dist[1]{d_{#1}}
\newcommand*\mcartes[1][n]{\nu_{#1}}
\newcommand*\marbres[1][n]{\mu_{#1}}
\newcommand*\mimage[1][,n]{\mu_{\bij #1}}

\begin{document}

\maketitle

\begin{abstract}
On d\'emontre que les quadrangulations al\'eatoires infinies
uniformes d\'efinies respectivement par Chassaing-Durhuus et par
Krikun ont la même loi.
\end{abstract}

\selectlanguage{english}
\begin{abstract}
We prove that the uniform infinite random quadrangulations defined
respectively by Chassaing-Durhuus and Krikun have the same distribution.
\end{abstract}
\selectlanguage{french}

\noindent {\bf Classification AMS:} 60C05, 60J80, 05C30.

\section{Introduction}

Les cartes planaires ont \'et\'e \'etudi\'ees dans un premier temps
par Tutte \cite{Tutte} dans une optique purement combinatoire. Ces
objets sont maintenant utilis\'es en physique th\'eorique  comme
mod\`eles de g\'eom\'etrie al\'eatoire dans la th\'eorie de la
gravit\'e quantique en deux dimensions (voir en particulier le livre \cite{2DQG}). Un
outil puissant pour l'\'etude de ces objets vient du codage des
cartes planaires par certaines classes d'arbres, introduit dans  un
premier temps par \cite{CV}, puis d\'evelopp\'e dans la th\`ese de
Schaeffer \cite{Schaeffer} (voir aussi Bouttier, Di Francesco, Guitter
\cite{BDG} pour une
g\'en\'eralisation de ce codage). Cette correspondance entre les
cartes planaires et les arbres permet de d\'ecrire les propri\'et\'es de
limites d'\'echelle de grandes cartes al\'eatoires gr\^ace \`a celles
des arbres continus
al\'eatoires (voir le travail de Chassaing et Schaeffer \cite{CS}), et
de d\'efinir une carte brownienne (voir Marckert et Mokkadem
\cite{MaMo}), carte al\'eatoire continue supposée \^etre la
limite d'\'echelle de nombreuses classes de cartes al\'eatoires (voir
les travaux de Marckert et Miermont \cite{MaMi},
Le Gall \cite{LG}, Le Gall et Paulin \cite{lGP}). Cette approche donne
de nouvelles propri\'et\'es asymptotiques des grandes cartes planaires.

\bigskip

Un autre point de vue consiste en l'étude de cartes planaires
infinies, plus particulièrement en l'étude de mesures de probabilité,
uniformes en un certain sens, sur des ensembles de cartes
infinies. Ceci a été fait par Angel et Schramm \cite{AS}, qui
introduisent une triangulation infinie uniforme du plan, dont les
propriétés sont étudiées par Angel \cite{A1,A2} et par Krikun \cite{Krikun}.

Chassaing et Durhuus construisent dans \cite{CD} un arbre
al\'eatoire bien \'etiquet\'e infini de loi uniforme comme la limite
locale d'arbres finis de loi uniforme parmi les arbres ayant une
taille donn\'ee. Ils \'etendent
ensuite la bijection de Schaeffer aux arbres et quadrangulations
infinis pour pouvoir d\'efinir \`a partir de cet arbre infini une
quadrangulation infinie al\'eatoire. La loi de
cette quadrangulation infinie est uniforme dans le sens o\`u elle est
la limite des probabilit\'es uniformes sur les quadrangulations ayant une
taille donn\'ee, pour la topologie image de la topologie naturelle
sur les arbres bien \'etiquet\'es par la bijection de Schaeffer. Cette
topologie sur les quadrangulations ne correspond cependant pas \`a la
topologie de la convergence locale, plus naturelle pour des objets
combinatoires comme les quadrangulations.

Krikun \cite{Kr} construit sans passer par les arbres une
quadrangulation infinie al\'eatoire comme la limite faible locale de
quadrangulations al\'eatoires de loi uniforme parmi les
quadrangulations d'une taille donn\'ee. Cette carte infinie al\'eatoire est
donc aussi en un sens naturel une quadrangulation infinie de loi
uniforme. Cependant la construction de Krikun est a priori
diff\'erente de celle de Chassaing et Durhuus.

\bigskip

Notre résultat principal, le théorème \ref{egalite}, montre que les deux d\'efinitions
donn\'ees respectivement par Chassaing-Durhuus et par Krikun sont
\'equivalentes. Cela justifie le nom de quadrangulation infinie
uniforme pour la carte al\'eatoire infinie obtenue par l'une ou
l'autre des d\'efinitions. Ce résultat est démontré dans
la partie \ref{sec:egalite} ci-dessous. La partie
\ref{sec:preliminaires} donne quelques préliminaires, concernant en
particulier la bijection de Schaeffer étendue entre arbres et
quadrangulations infinis. La partie \ref{sec:quadrangulations} rapelle les deux
définitions de la quadrangulation infinie uniforme.

\section{Pr\'eliminaires}
\label{sec:preliminaires}

\subsection{Arbres bien \'etiquet\'es}

Pour d\'ecrire les arbres planaires nous utilisons le formalisme de
Neveu introduit dans \cite{Neveu}. Soit
\[ \mathcal{U} = \bigcup_{n = 0}^{\infty} \mathbb{N}^n \]
o\`u $\mathbb{N}=\{1,2,\ldots \}$ et $\mathbb{N}^0 = \{ \emptyset \}$. La \emph{g\'en\'eration}
de $u = (u_1,\ldots, u_n) \in \mathbb{N}^n$ est $|u|=n$. Si $u = (u_1, \ldots , u_m)$ et
$v = (v_1, \ldots , v_n)$ sont deux \'el\'ements de $\mathcal{U}$, $uv = (u_1, \ldots,u_m, v_1, \ldots , v_n)$ est la concat\'enation de $u$ et $v$ (en particulier $\emptyset
u = u \emptyset = u$). Si $v$ est de la forme $uj$ avec $u \in \mathcal{U}$ et $j \in
\mathbb{N}$, on dit que $u$ est le \emph{p\`ere} de $v$ et aussi que $v$ est
un \emph{enfant} de $u$. De mani\`ere plus g\'en\'erale, si $v$ est de la
forme $uw$ pour $u, \, w \in \mathcal{U}$, on dit que $u$ est un
\emph{anc\^etre} de $v$ et aussi que $v$ est un \emph{descendant} de $u$.

\bigskip

Un \emph{arbre planaire enracin\'e} (ou plus simplement par la suite un
arbre) $\tau$ est un sous ensemble de $\mathcal{U}$ tel que:
\begin{enumerate}
\item $\emptyset \in \tau$ (on dit que $\emptyset$ est la \emph{racine} de $\tau$);
\item si $v \in \tau$ et $v \neq \emptyset$, le p\`ere de $v$ appartient \`a $\tau$;
\item pour tout $u \in \tau$ il existe un entier $k_u(\tau) \geqslant
  0$ tel que $uj \in \tau$ si et seulement si $1 \leqslant j \leqslant k_u(\tau)$.
\end{enumerate}
Les arêtes de $\tau$ sont les couples $(u,v)$ o\`u $u,v \in \tau$ et
$u$ est le père de $v$.
$|\tau|$ d\'esigne le nombre d'ar\^etes de $\tau$ et est appel\'e taille de
$\tau$; $h(\tau)$ d\'esigne la g\'en\'eration maximale des sommets de $\tau$
et est appel\'e hauteur de $\tau$.
On note $\mathbf{T}_n$ l'ensemble des arbres planaires enracin\'es de
taille $n$, $\mathbf{T}_{\infty}$ l'ensemble des arbres infinis. $\mathbf{T}
= \bigcup_{n = 0}^{\infty} \mathbf{T}_n$ est l'ensemble des arbres
planaires enracin\'es finis; $\overline{\mathbf{T}} =
\mathbf{T} \cup \mathbf{T}_{\infty}$ est l'ensemble des arbres planaires
enracin\'es, finis ou infinis.

\bigskip

Un \emph{arbre planaire enracin\'e \'etiquet\'e} (ou arbre spatial, ou plus
simplement arbre \'etiquet\'e) est une paire $\omega =
(\tau, (\ell(u))_{u \in \tau})$ constitu\'ee d'un arbre planaire enracin\'e
$\tau$ et d'une collection d'\'etiquettes attach\'ees aux sommets de $\tau$
telles que $\ell(u) \in \mathbb{Z}$ pour tout $u \in \tau$ et
$|\ell(u) - \ell(v)| \leqslant 1$ si $(u,v)$ est une arête de $\tau$.
Par convention, $|\omega|= |\tau|$ est la taille de $\omega$. 
Pour $l \in \mathbb{Z}$ on note
$\mathbf{T}_n^{(l)}$ l'ensemble des arbres planaires enracin\'es et
\'etiquet\'es \`a $n$ ar\^etes tels que $\ell(\emptyset) = l$ ainsi
que $\mathbf{T}^{(l)}$, $\mathbf{T}_{\infty}^{(l)}$ et
$\overline{\mathbf{T}}^{(l)}$ les ensembles de tels arbres
respectivement finis, infinis et de taille quelconque.

Si de plus on a $\ell(u) \geqslant l$ pour tout sommet $u$ de $\tau$, on
dit que $\omega$ est \emph{$l$-bien \'etiquet\'e}.
Les ensembles d'arbres $l$-bien \'etiquet\'es correspondants seront not\'es
$\arbre[n]^{(l)}$, $\arbre[{}]^{(l)}$, $\arbre^{(l)}$ et
$\arbres^{(l)}$. Pour $l = 1$ on parle d'\emph{arbres bien \'etiquet\'es} et
on notera les ensembles d'arbres bien \'etiquet\'es plus simplement $\arbre[n]$,
$\arbre[{}]$, $\arbre$ et $\arbres$.

\bigskip

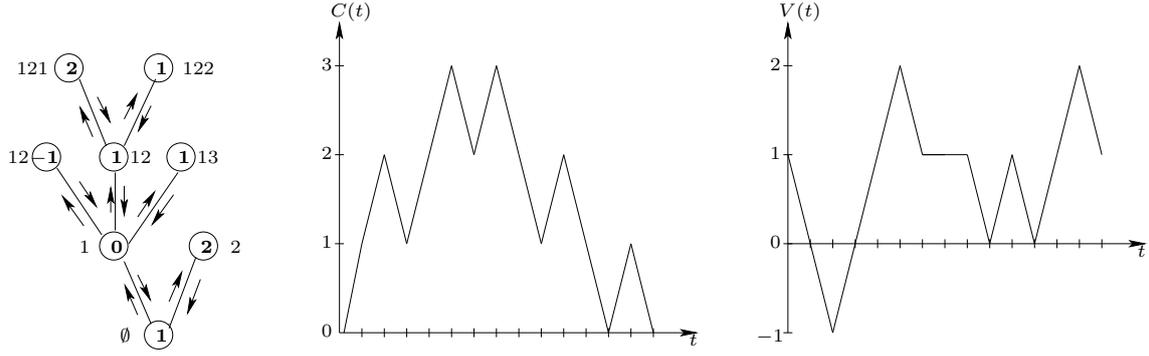
\begin{figure}[!t]
\begin{center}
\begin{picture}(0,0)%
\includegraphics{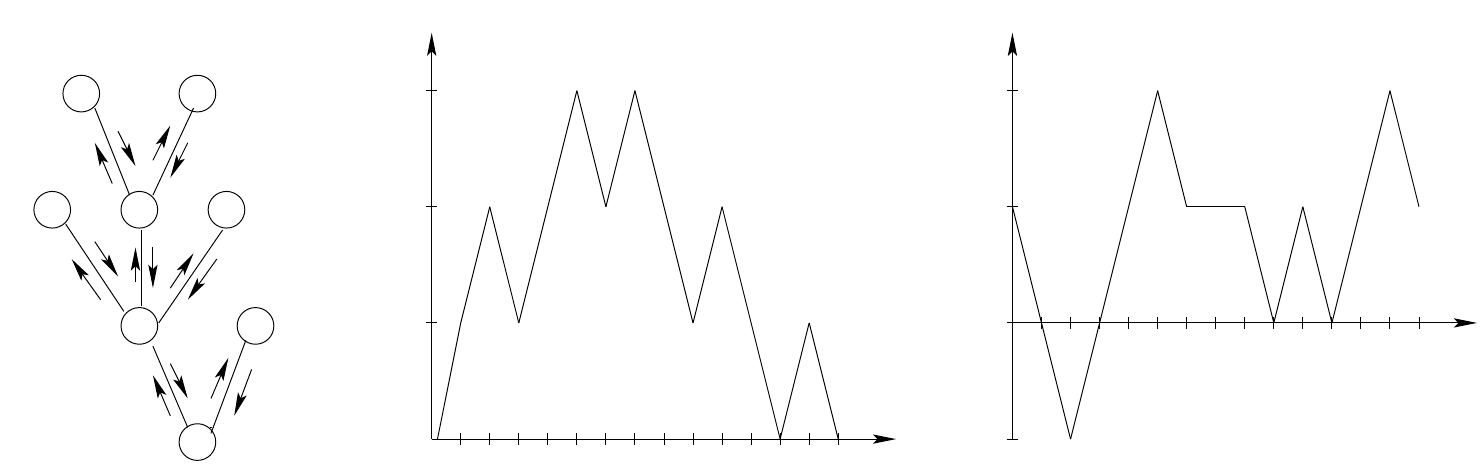}%
\end{picture}%
\setlength{\unitlength}{2445sp}%
\begingroup\makeatletter\ifx\SetFigFontNFSS\undefined%
\gdef\SetFigFontNFSS#1#2#3#4#5{%
  \reset@font\fontsize{#1}{#2pt}%
  \fontfamily{#3}\fontseries{#4}\fontshape{#5}%
  \selectfont}%
\fi\endgroup%
\begin{picture}(11457,3582)(3856,-6340)
\put(11611,-2941){\makebox(0,0)[lb]{\smash{{\SetFigFontNFSS{7}{8.4}{\rmdefault}{\mddefault}{\updefault}{\color[rgb]{0,0,0}$V(t)$}%
}}}}
\put(7111,-2941){\makebox(0,0)[lb]{\smash{{\SetFigFontNFSS{7}{8.4}{\rmdefault}{\mddefault}{\updefault}{\color[rgb]{0,0,0}$C(t)$}%
}}}}
\put(7021,-5281){\makebox(0,0)[lb]{\smash{{\SetFigFontNFSS{7}{8.4}{\rmdefault}{\mddefault}{\updefault}{\color[rgb]{0,0,0}$1$}%
}}}}
\put(7021,-6181){\makebox(0,0)[lb]{\smash{{\SetFigFontNFSS{7}{8.4}{\rmdefault}{\mddefault}{\updefault}{\color[rgb]{0,0,0}$0$}%
}}}}
\put(10711,-6271){\makebox(0,0)[lb]{\smash{{\SetFigFontNFSS{7}{8.4}{\rmdefault}{\mddefault}{\updefault}{\color[rgb]{0,0,0}$t$}%
}}}}
\put(7021,-3481){\makebox(0,0)[lb]{\smash{{\SetFigFontNFSS{7}{8.4}{\rmdefault}{\mddefault}{\updefault}{\color[rgb]{0,0,0}$3$}%
}}}}
\put(7021,-4381){\makebox(0,0)[lb]{\smash{{\SetFigFontNFSS{7}{8.4}{\rmdefault}{\mddefault}{\updefault}{\color[rgb]{0,0,0}$2$}%
}}}}
\put(3961,-3526){\makebox(0,0)[lb]{\smash{{\SetFigFontNFSS{7}{8.4}{\rmdefault}{\mddefault}{\updefault}{\color[rgb]{0,0,0}$121$}%
}}}}
\put(5626,-3526){\makebox(0,0)[lb]{\smash{{\SetFigFontNFSS{7}{8.4}{\rmdefault}{\mddefault}{\updefault}{\color[rgb]{0,0,0}$122$}%
}}}}
\put(3871,-4426){\makebox(0,0)[lb]{\smash{{\SetFigFontNFSS{7}{8.4}{\rmdefault}{\mddefault}{\updefault}{\color[rgb]{0,0,0}$12$}%
}}}}
\put(4591,-5326){\makebox(0,0)[lb]{\smash{{\SetFigFontNFSS{7}{8.4}{\rmdefault}{\mddefault}{\updefault}{\color[rgb]{0,0,0}$1$}%
}}}}
\put(5806,-5326){\makebox(0,0)[lb]{\smash{{\SetFigFontNFSS{7}{8.4}{\rmdefault}{\mddefault}{\updefault}{\color[rgb]{0,0,0}$\mathbf{2}$ \, $2$}%
}}}}
\put(5356,-6226){\makebox(0,0)[lb]{\smash{{\SetFigFontNFSS{7}{8.4}{\rmdefault}{\mddefault}{\updefault}{\color[rgb]{0,0,0}$\mathbf{1}$}%
}}}}
\put(15211,-5371){\makebox(0,0)[lb]{\smash{{\SetFigFontNFSS{7}{8.4}{\rmdefault}{\mddefault}{\updefault}{\color[rgb]{0,0,0}$t$}%
}}}}
\put(11521,-5281){\makebox(0,0)[lb]{\smash{{\SetFigFontNFSS{7}{8.4}{\rmdefault}{\mddefault}{\updefault}{\color[rgb]{0,0,0}$0$}%
}}}}
\put(11521,-4381){\makebox(0,0)[lb]{\smash{{\SetFigFontNFSS{7}{8.4}{\rmdefault}{\mddefault}{\updefault}{\color[rgb]{0,0,0}$1$}%
}}}}
\put(11521,-3481){\makebox(0,0)[lb]{\smash{{\SetFigFontNFSS{7}{8.4}{\rmdefault}{\mddefault}{\updefault}{\color[rgb]{0,0,0}$2$}%
}}}}
\put(4996,-6226){\makebox(0,0)[lb]{\smash{{\SetFigFontNFSS{7}{8.4}{\rmdefault}{\mddefault}{\updefault}{\color[rgb]{0,0,0}$\emptyset$}%
}}}}
\put(4906,-5326){\makebox(0,0)[lb]{\smash{{\SetFigFontNFSS{7}{8.4}{\rmdefault}{\mddefault}{\updefault}{\color[rgb]{0,0,0}$\mathbf{0}$}%
}}}}
\put(4906,-4426){\makebox(0,0)[lb]{\smash{{\SetFigFontNFSS{7}{8.4}{\rmdefault}{\mddefault}{\updefault}{\color[rgb]{0,0,0}$\mathbf{1}$ $12$}%
}}}}
\put(5581,-4426){\makebox(0,0)[lb]{\smash{{\SetFigFontNFSS{7}{8.4}{\rmdefault}{\mddefault}{\updefault}{\color[rgb]{0,0,0}$\mathbf{1}$ $13$}%
}}}}
\put(4456,-3526){\makebox(0,0)[lb]{\smash{{\SetFigFontNFSS{7}{8.4}{\rmdefault}{\mddefault}{\updefault}{\color[rgb]{0,0,0}$\mathbf{2}$}%
}}}}
\put(5356,-3526){\makebox(0,0)[lb]{\smash{{\SetFigFontNFSS{7}{8.4}{\rmdefault}{\mddefault}{\updefault}{\color[rgb]{0,0,0}$\mathbf{1}$}%
}}}}
\put(11386,-6226){\makebox(0,0)[lb]{\smash{{\SetFigFontNFSS{7}{8.4}{\rmdefault}{\mddefault}{\updefault}{\color[rgb]{0,0,0}$-1$}%
}}}}
\put(4096,-4426){\makebox(0,0)[lb]{\smash{{\SetFigFontNFSS{7}{8.4}{\rmdefault}{\mddefault}{\updefault}{\color[rgb]{0,0,0}$\mathbf{-1}$}%
}}}}
\end{picture}%
\end{center}
\caption{Un arbre \'etiquet\'e et ses fonctions de contour $(C,V)$.}
\label{figcontour}
\end{figure}

Un arbre planaire enracin\'e \'etiquet\'e $\omega = (\tau, \ell)$ peut
\^etre cod\'e par une paire de fonctions $(C,V)$, o\`u $C = (C(t), 0
\leqslant t \leqslant 2 |\omega|)$ est la fonction de \emph{contour} de
$\tau$ et $V = (V(t), 0 \leqslant t \leqslant 2 |\omega| )$ est
la fonction de \emph{contour spatial} de $\omega$. Pour d\'efinir la fonction
de contour de $\tau$ consid\'erons une particule qui, partant de la
racine, parcourt l'arbre le long de ses ar\^etes. Lorsqu'elle quitte un
sommet, la particule se d\'eplace vers le premier enfant non d\'ej\`a
visit\'e de ce sommet si c'est possible,
et retourne vers la racine dans le cas contraire.
On convient que chaque arête a pour longueur $1$ et que la particule
se déplace à vitesse unité. La particule met alors un temps \'egal \`a
$2|\tau|$ pour parcourir l'arbre en entier (chaque ar\^ete est
parcourue deux fois). Pour tout $t \in [0,2|\tau| ]$, $C(t)$ est la
distance de la particule \`a la racine. La définition de la fonction
de contour devrait être claire sur la figure \ref{figcontour}.
D\'efinissons maintenant la
fonction de contour spatial.  Si $t \in [0,2|\tau|]$ est
un entier, la particule est \`a l'instant $t$ sur un sommet $v$ et on pose $V(t) =
\ell(v)$, puis on compl\`ete la d\'efinition de $V$ par une interpolation
lin\'eaire entre deux entiers successifs (voir la figure \ref{figcontour})
Un arbre \'etiquet\'e $\omega$ est clairement
d\'etermin\'e de fa\c{c}on unique par ses fonctions de contour.

\bigskip

Nous aurons besoin des d\'efinitions suivantes emprunt\'ees \`a \cite{CS}:

\begin{definition}
Une \emph{colonne} d'un arbre bien \'etiquet\'e infini $\omega$ est une suite infinie
$(s_n)_{n \geqslant 0}$ de sommets qui commence au sommet racine
et telle que pour tout $n > 0$, $s_n$ est un fils de
$s_{n-1}$. De mani\`ere \'equivalente, une colonne de $\omega$ est un
sous-arbre lin\'eaire infini de $\omega$ ayant le m\^eme
sommet racine que $\omega$.

On note $\spine$ le sous ensemble de $\arbre$ constitu\'e des arbres
ayant une unique colonne.
\end{definition}

Pour tout $l \in \mathbb{N}$, on note $N_l$ la fonction qui \`a un arbre bien
\'etiquet\'e $\omega$ associe le nombre de sommets de $\omega$ ayant pour
\'etiquette $l$. On d\'efinit alors le sous-ensemble de $\spine$
suivant: 

\begin{definition}
\label{colonneens}
On pose
\[ \spinelab = \left\{ \omega \in \spine : \, \forall l \geqslant 1 , \, N_l(\omega)
    < \infty \right\} \cup \mathbb{T} ,\]
c'est l'ensemble des arbres bien \'etiquet\'es, finis ou infinis, ayant une
unique colonne dans le cas o\`u il sont infinis, et dont toutes les
\'etiquettes n'apparaissent qu'un nombre fini de fois.
\end{definition}

Introduisons finalement les notations suivantes, qui nous seront utiles
par la suite; si $\omega = (\tau, \ell)$ est un arbre \'etiquet\'e,
$|\omega| =|\tau|$ et $h(\omega) = h(\tau)$ d\'esignent la taille et
la hauteur de $\omega$, et, pour tout entier $S \geqslant 0$,
$g_S(\omega)$ est l'ensemble des sommets
de g\'en\'eration $S$ de $\omega$.

\subsection{Quadrangulations}

Une \emph{carte planaire} est le plongement propre dans la sph\`ere orient\'ee $\mathbb{S}^2$
d'un graphe planaire connexe et fini. Ce plongement est
consid\'er\'e \`a hom\'eomorphisme de la sph\`ere conservant
l'orientation pr\`es. Les composantes connexes du compl\'ementaire du
graphe sont appel\'ees faces de la carte.
Une carte planaire est dite \emph{enracin\'ee} si une
ar\^ete orient\'ee, appell\'ee \emph{racine}, est sp\'ecifi\'ee. Si de plus
toutes les faces de la carte ont quatre ar\^etes adjacentes (une
ar\^ete totalement inscrite dans une face est compt\'ee deux fois), la
carte est une \emph{quadrangulation} (enracin\'ee).

On note $\carte[n]$ l'ensemble des quadrangulations enracin\'ees \`a $n$
faces et $\carte[{}] = \bigcup_{n = 0}^{\infty} \carte[n]$
l'ensemble des quadrangulations finies. Chaque ensemble $\carte[n]$
est en correspondance bijective avec l'ensemble $\arbre[n]$ gr\^ace \`a la
bijection de Schaeffer \cite{CV,Schaeffer}. Nous allons voir par la suite
comment \'etendre cette correspondance pour transformer un arbre de
$\spinelab$ en une quadrangulation enracin\'ee infinie de $\mathbb{S}^2$. Cependant
aucune bijection n'est \'etablie entre les arbres bien \'etiquet\'es
infinis et les quadrangulations enracin\'ees infinies.

\bigskip

Avant de voir comment associer une quadrangulation \`a un arbre de $\spinelab$, nous
devons d\'efinir plus pr\'ecis\'ement la notion de quadrangulation
enracin\'ee infinie. Cette d\'efinition nous conduira \`a introduire
quelques propri\'et\'es de r\'egularit\'e. Nous reprenons dans ce but les
d\'efinitions faites dans \cite{AS,CD}.

\begin{definition}
Une \emph{carte planaire infinie} $\mathcal{M}$ est le plongement dans la
sph\`ere $\mathbb{S}^2$ d'un graphe infini $G$, connexe et de sommets de
degr\'es finis, tel que les ar\^etes sont repr\'esent\'ees par des portions
de courbes r\'eguli\`eres qui ne s'intersectent entre elles que sur les
sommets qu'elles ont en commun.
\end{definition}

Avec cette d\'efinition, l'arbre lin\'eaire infini plong\'e dans
$\mathbb{S}^2$ de sorte qu'il occupe un cercle entier est une carte
planaire. Nous ne souhaitons pas accepter ce plongement comme valide car il
donne une carte planaire \`a deux <<faces>> l\`a o\`u nous en voulons une
seule. Ceci motive la d\'efinition suivante:

\begin{definition}
\label{acceptable}
Une carte planaire infinie \emph{acceptable} $\mathcal{M}$ est une carte
planaire infinie v\'erifiant la propri\'et\'e de r\'egularit\'e
suivante:

Si $(p_i)_{i \in \mathbb{N}}$ est une suite de points de $\mathcal{M}$
consid\'er\'ee comme l'union de ses ar\^etes dans $\mathbb{S}^2$, et telle
que pour $i \neq j$, $p_i$ et $p_j$ sont sur deux ar\^etes
diff\'erentes, alors la suite n'a aucun point d'accumulation dans
$\mathcal{M}$.
\end{definition}

Cette propri\'et\'e de r\'egularit\'e va nous permettre de d\'efinir ce que
sont les faces d'une carte planaire infinie acceptable. Soit en effet $C$
une courbe ferm\'ee et continue d'une carte planaire acceptable
$\mathcal{M}$, compos\'ee d'une suite d'ar\^etes. La propri\'et\'e de
r\'egularit\'e assure que le nombre d'ar\^etes qui constituent $C$ est
fini, ainsi le compl\'ementaire de $C$ dans $\mathbb{S}^2$ se d\'ecompose en un
nombre fini de composantes connexes. Si une de ces composantes connexes
contient un nombre fini de sommets de $\mathcal{M}$, alors la partie de
$\mathcal{M}$ \`a l'int\'erieur de cette composante est une
carte planaire finie. Les faces de cette carte planaire finie sont
appel\'ees faces de $\mathcal{M}$. Par d\'efinition, les faces de
$\mathcal{M}$ sont toutes obtenues de cette mani\`ere.
Elles sont en particulier toutes born\'ees par une boucle polygonale
compos\'ee d'un nombre fini de sommets.

Avec cette d\'efinition, toutes les ar\^etes d'une carte planaire infinie
acceptable ne sont pas n\'ecessairement adjacentes \`a une face. C'est le
cas par exemple des arbres infinis, qui ont une seule << face >> de degr\'e
infini (laquelle n'est pas une face au sens de la définition
précédente). Nous aurons donc aussi besoin de la d\'efinition suivante:

\begin{definition}
Une carte planaire infinie \emph{r\'eguli\`ere} est une carte planaire infinie
acceptable telle que toute ar\^ete de la carte est partag\'ee par
exactement deux faces, ou bien appara\^it deux fois dans la fronti\`ere d'une face.
\end{definition}

Comme dans le cas des cartes planaires finies, deux cartes planaires
infinies régulières sont identifiées si l'on passe de l'une à l'autre
par un homéomorphisme préservant l'orientation.

\begin{remarque}
Avec les d\'efinitions que nous avons donn\'ees, un arbre infini peut se
plonger dans $\mathbb{S}^2$ en une carte planaire infinie acceptable mais
pas en une carte planaire infinie r\'eguli\`ere.
\end{remarque}

\begin{definition}
\label{Qinf}
Une \emph{quadrangulation infinie} est une carte planaire infinie r\'eguli\`ere
dont toutes les faces sont d\'elimit\'ees par des polygones \`a quatre
cot\'es. Une quadrangulation enracin\'ee est une quadrangulation munie d'une
ar\^ete orient\'ee $(i_0,i_1)$ distingu\'ee, appel\'ee racine de la
quadrangulation; $i_0$ \'etant appel\'e sommet racine.
On note $\cartes$ l'ensemble des quadrangulations enracin\'ees finies ou
infinies, avec la
d\'ecomposition \'evidente $\cartes = \carte[{}] \cup \carte$
\end{definition}

\subsection{Bijection de Schaeffer \'etendue}

\label{sch}

Nous allons maintenant pr\'esenter la bijection de Schaeffer
\'etendue aux arbres infinis de l'ensemble $\spinelab$.
Nous renvoyons \`a \cite{CD} (section 6) pour plus de d\'etails et pour la preuve des
r\'esultats qui suivent.

Pour construire l'image $\bij(\omega)$ d'un arbre infini $\omega \in
\spinelab$, nous avons besoin de propri\'et\'es de r\'egularit\'e et
surtout de l'unicit\'e de son plongement dans $\mathbb{S}^2$. Un arbre infini
peut en effet \^etre plong\'e dans la sph\`ere de plusieurs mani\`eres non
hom\'eomorphes. Cependant, il n'existe \`a hom\'eomorphisme pr\`es
qu'une seule fa\c{c}on de plonger un
arbre de $\spinelab$ de mani\`ere \`a obtenir une carte planaire
acceptable pour laquelle il existe $c \in \mathbb{S}^2$ tel que, quelle que 
soit la suite de points $(p_i)_{i \in \mathbb{N}}$ v\'erifiant les
conditions de la d\'efinition \ref{acceptable} et ayant effectivement
un point d'accumulation dans $\mathbb{S}^2$, ce point est $c$.

Soit donc $\omega \in \spinelab$, et fixons un plongement de $\omega$
dans $\mathbb{S}^2$ comme carte planaire acceptable.
Soit $c$ le seul point d'accumulation possible
des suites de la d\'efinition \ref{acceptable}. On note $F_0$ le
compl\'ementaire de $\omega$ dans $\mathbb{S}^2$.

\begin{definition}
Un \emph{coin} de $F_0$  est le secteur d\'efini autour d'un sommet de $\omega$
par deux ar\^etes cons\'ecutives. L'\'etiquette d'un coin est l'\'etiquette
du sommet qui lui correspond.
\end{definition}

Il est facile de voir qu'un sommet de degr\'e $d$ d\'efinit $d$ coins, et
que $\omega$ a un nombre fini $C_l(\omega) \geqslant N_l(\omega)$ de coins
d'\'etiquette $l$. Nous allons construire $\bij(\omega)$ en trois \'etapes.

\bigskip

\begin{figure}[!t]
\begin{center}
\begin{picture}(0,0)%
\includegraphics{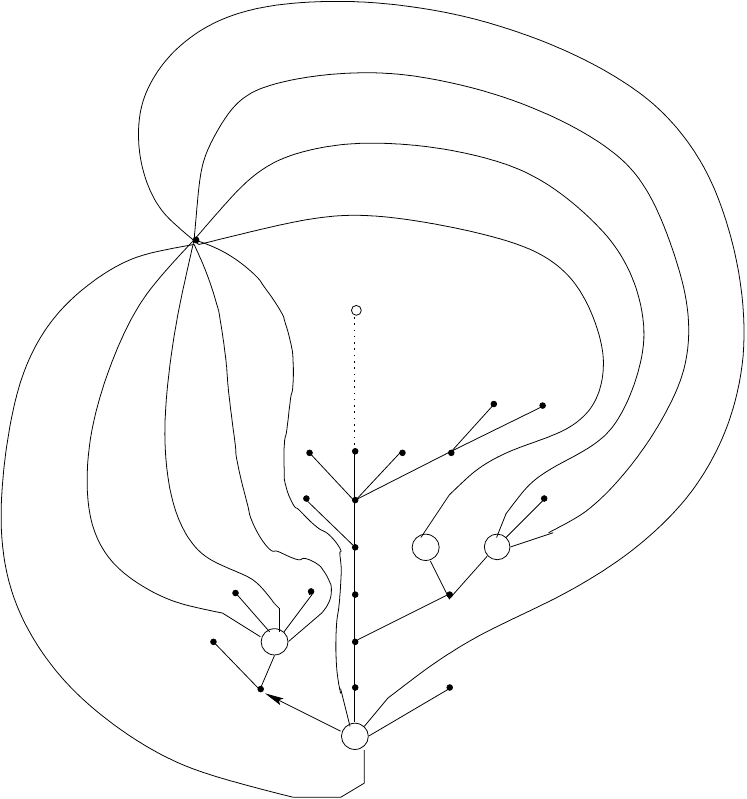}%
\end{picture}%
\setlength{\unitlength}{1989sp}%
\begingroup\makeatletter\ifx\SetFigFontNFSS\undefined%
\gdef\SetFigFontNFSS#1#2#3#4#5{%
  \reset@font\fontsize{#1}{#2pt}%
  \fontfamily{#3}\fontseries{#4}\fontshape{#5}%
  \selectfont}%
\fi\endgroup%
\begin{picture}(7099,7604)(1121,-7858)
\put(4591,-3166){\makebox(0,0)[lb]{\smash{{\SetFigFontNFSS{6}{7.2}{\familydefault}{\mddefault}{\updefault}{\color[rgb]{0,0,0}$c$}%
}}}}
\put(2701,-2356){\makebox(0,0)[lb]{\smash{{\SetFigFontNFSS{6}{7.2}{\rmdefault}{\mddefault}{\updefault}{\color[rgb]{0,0,0}$v_0$}%
}}}}
\put(5356,-3301){\makebox(0,0)[lb]{\smash{{\SetFigFontNFSS{6}{7.2}{\rmdefault}{\mddefault}{\updefault}{\color[rgb]{0,0,0}$F_{\infty}$}%
}}}}
\put(7021,-7216){\makebox(0,0)[lb]{\smash{{\SetFigFontNFSS{6}{7.2}{\familydefault}{\mddefault}{\updefault}{\color[rgb]{0,0,0}$F_8$}%
}}}}
\put(3646,-691){\makebox(0,0)[lb]{\smash{{\SetFigFontNFSS{6}{7.2}{\familydefault}{\mddefault}{\updefault}{\color[rgb]{0,0,0}$F_1$}%
}}}}
\put(4006,-1321){\makebox(0,0)[lb]{\smash{{\SetFigFontNFSS{6}{7.2}{\familydefault}{\mddefault}{\updefault}{\color[rgb]{0,0,0}$F_2$}%
}}}}
\put(4501,-1951){\makebox(0,0)[lb]{\smash{{\SetFigFontNFSS{6}{7.2}{\familydefault}{\mddefault}{\updefault}{\color[rgb]{0,0,0}$F_3$}%
}}}}
\put(3511,-3481){\makebox(0,0)[lb]{\smash{{\SetFigFontNFSS{6}{7.2}{\familydefault}{\mddefault}{\updefault}{\color[rgb]{0,0,0}$F_4$}%
}}}}
\put(3016,-4426){\makebox(0,0)[lb]{\smash{{\SetFigFontNFSS{6}{7.2}{\familydefault}{\mddefault}{\updefault}{\color[rgb]{0,0,0}$F_5$}%
}}}}
\put(2296,-4786){\makebox(0,0)[lb]{\smash{{\SetFigFontNFSS{6}{7.2}{\familydefault}{\mddefault}{\updefault}{\color[rgb]{0,0,0}$F_6$}%
}}}}
\put(1846,-5956){\makebox(0,0)[lb]{\smash{{\SetFigFontNFSS{6}{7.2}{\familydefault}{\mddefault}{\updefault}{\color[rgb]{0,0,0}$F_7$}%
}}}}
\put(4456,-7351){\makebox(0,0)[lb]{\smash{{\SetFigFontNFSS{6}{7.2}{\familydefault}{\mddefault}{\updefault}{\color[rgb]{0,0,0}$\mathbf{1}$}%
}}}}
\put(3691,-6451){\makebox(0,0)[lb]{\smash{{\SetFigFontNFSS{6}{7.2}{\familydefault}{\mddefault}{\updefault}{\color[rgb]{0,0,0}$\mathbf{1}$}%
}}}}
\put(5131,-5551){\makebox(0,0)[lb]{\smash{{\SetFigFontNFSS{6}{7.2}{\familydefault}{\mddefault}{\updefault}{\color[rgb]{0,0,0}$\mathbf{1}$}%
}}}}
\put(5806,-5551){\makebox(0,0)[lb]{\smash{{\SetFigFontNFSS{6}{7.2}{\familydefault}{\mddefault}{\updefault}{\color[rgb]{0,0,0}$\mathbf{1}$}%
}}}}
\end{picture}%
\end{center}
\caption{\'Etape 1: on relie $v_0$ aux coins d'\'etiquette $1$.}
\label{coins1}
\end{figure}
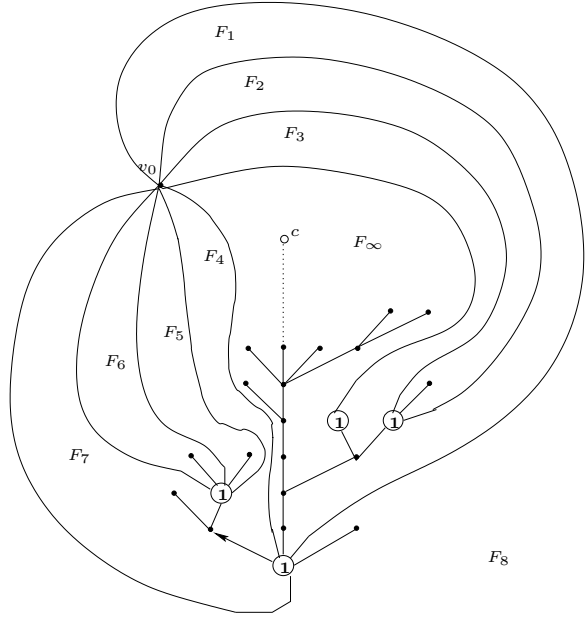

\begin{etape}[voir figure \ref{coins1}] 
On ajoute un sommet $v_0$ d'\'etiquette $0$ dans $F_0
\setminus \{ c \}$ que l'on relie \`a chacun des $C_1(\omega)$ coins
d'\'etiquette $1$. La racine de la carte
planaire acceptable $\mathcal{M}_0$ obtenue est l'ar\^ete joignant $v_0$ au
coin pr\'ec\'edant la racine de $\omega$.
\end{etape}

\begin{remarque}
Il est possible de relier $v_0$ aux coins d'\'etiquette $1$ car $\omega$ a
exactement une colonne.
\end{remarque}

Apr\`es cette 1\iere\ \'etape, une carte planaire acceptable enracin\'ee
$\mathcal{M}_0$ poss\'edant $C_1(\omega) - 1$ faces est obtenue (au sens des
d\'efinitions donn\'ees pr\'ec\'edemment, en particulier les faces
sont finies). Remarquons que
chacune des faces de $\mathcal{M}_0$ poss\`ede un unique coin d'\'etiquette
$0$ et deux coins d'\'etiquette $1$ qui sont du m\^eme c\^ot\'e de la colonne
de $\omega$. Ces faces sont d\'elimit\'ees chacune par les deux ar\^etes joignant
$v_0$ aux deux coins d'\'etiquette $1$ ainsi que, dans le cas o\`u ces coins
d'\'etiquette $1$ appartiennent \`a deux sommets distincts, les ar\^etes
figurant sur l'unique chemin injectif joignant ces deux
sommets d'\'etiquette $1$ dans l'arbre.

Il est naturel de consid\'erer le compl\'ementaire de $\mathcal{M}_0$ et de ses faces
comme une face suppl\'ementaire de degr\'e infini. On apellera cette face
$F_{\infty}$. Elle poss\`ede un unique coin d'\'etiquette $0$ et
deux coins d'\'etiquette $1$, qui sont de part et d'autre de la colonne
de $\omega$. Ces deux coins sont de plus les derniers coins
d'\'etiquette $1$ visit\'es lors d'un contour des c\^ot\'es
gauche et droit de $\omega$. $F_{\infty}$ est donc d\'elimit\'ee par les
deux ar\^etes joignant $v_0$ \`a ces coins d'\'etiquette $1$ et les ar\^etes
figurant sur l'unique chemin injectif joignant les deux sommets
correspondants dans l'arbre. La colonne de $\omega$ est, \`a l'exception d'un
nombre fini d'ar\^etes, \`a l'int\'erieur de cette face.

\bigskip

\begin{figure}[!t]
\begin{center}
\begin{picture}(0,0)%
\includegraphics{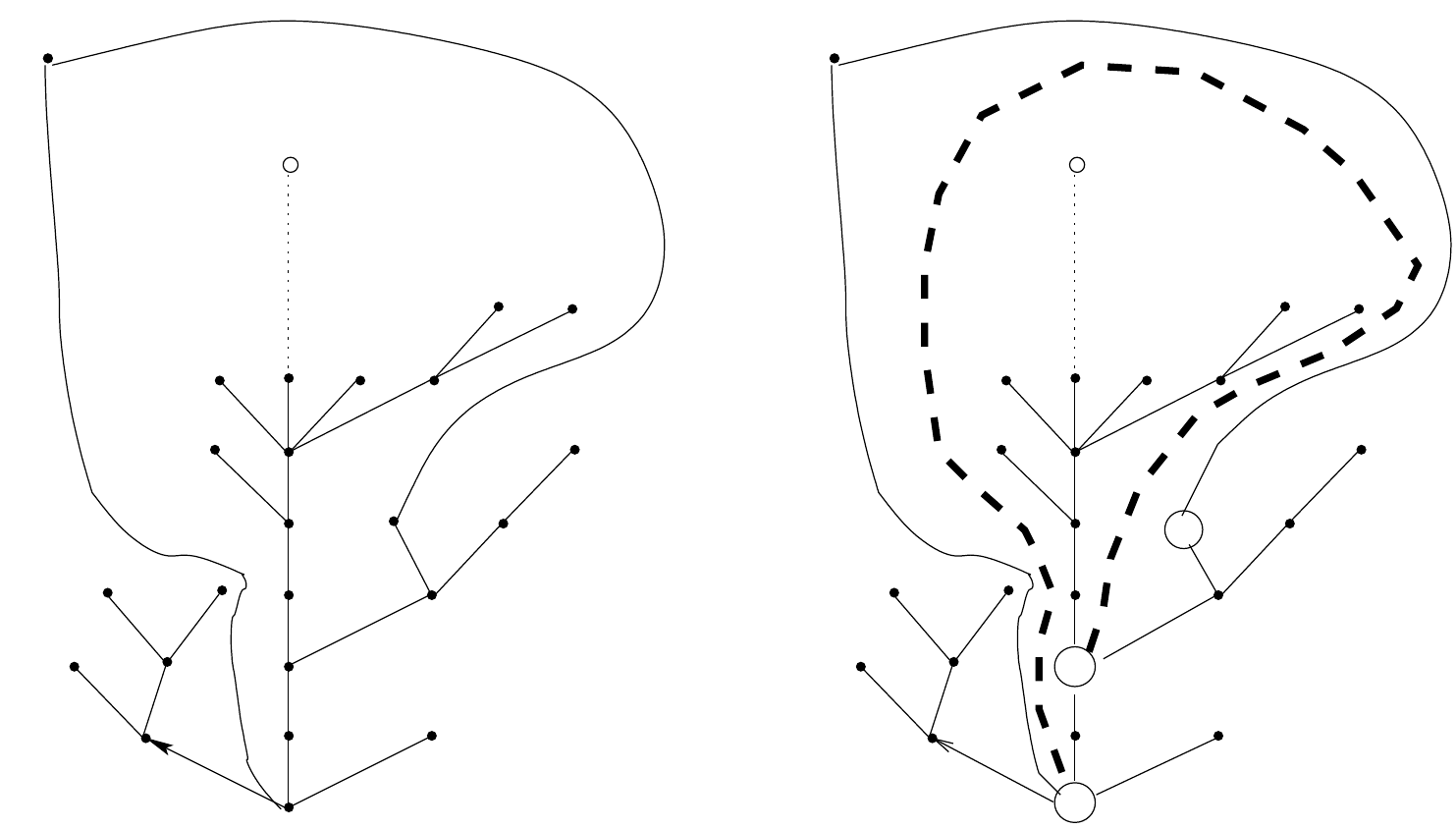}%
\end{picture}%
\setlength{\unitlength}{3066sp}%
\begingroup\makeatletter\ifx\SetFigFontNFSS\undefined%
\gdef\SetFigFontNFSS#1#2#3#4#5{%
  \reset@font\fontsize{#1}{#2pt}%
  \fontfamily{#3}\fontseries{#4}\fontshape{#5}%
  \selectfont}%
\fi\endgroup%
\begin{picture}(9145,5211)(2686,-7384)
\put(4591,-3166){\makebox(0,0)[lb]{\smash{{\SetFigFontNFSS{9}{10.8}{\familydefault}{\mddefault}{\updefault}{\color[rgb]{0,0,0}$c$}%
}}}}
\put(2701,-2356){\makebox(0,0)[lb]{\smash{{\SetFigFontNFSS{9}{10.8}{\familydefault}{\mddefault}{\updefault}{\color[rgb]{0,0,0}$v_0$}%
}}}}
\put(5356,-3301){\makebox(0,0)[lb]{\smash{{\SetFigFontNFSS{9}{10.8}{\familydefault}{\mddefault}{\updefault}{\color[rgb]{0,0,0}$F_{\infty}$}%
}}}}
\put(4591,-5956){\makebox(0,0)[lb]{\smash{{\SetFigFontNFSS{9}{10.8}{\familydefault}{\mddefault}{\updefault}{\color[rgb]{0,0,0}$4$}%
}}}}
\put(4591,-5596){\makebox(0,0)[lb]{\smash{{\SetFigFontNFSS{9}{10.8}{\familydefault}{\mddefault}{\updefault}{\color[rgb]{0,0,0}$5$}%
}}}}
\put(4636,-5146){\makebox(0,0)[lb]{\smash{{\SetFigFontNFSS{9}{10.8}{\familydefault}{\mddefault}{\updefault}{\color[rgb]{0,0,0}$6$}%
}}}}
\put(5401,-4741){\makebox(0,0)[lb]{\smash{{\SetFigFontNFSS{9}{10.8}{\familydefault}{\mddefault}{\updefault}{\color[rgb]{0,0,0}$7$}%
}}}}
\put(6436,-4066){\makebox(0,0)[lb]{\smash{{\SetFigFontNFSS{9}{10.8}{\familydefault}{\mddefault}{\updefault}{\color[rgb]{0,0,0}$8$}%
}}}}
\put(5716,-4021){\makebox(0,0)[lb]{\smash{{\SetFigFontNFSS{9}{10.8}{\familydefault}{\mddefault}{\updefault}{\color[rgb]{0,0,0}$10$}%
}}}}
\put(5221,-5866){\makebox(0,0)[lb]{\smash{{\SetFigFontNFSS{9}{10.8}{\familydefault}{\mddefault}{\updefault}{\color[rgb]{0,0,0}$2$}%
}}}}
\put(3871,-4426){\makebox(0,0)[lb]{\smash{{\SetFigFontNFSS{9}{10.8}{\familydefault}{\mddefault}{\updefault}{\color[rgb]{0,0,0}$-9$}%
}}}}
\put(3691,-4921){\makebox(0,0)[lb]{\smash{{\SetFigFontNFSS{9}{10.8}{\familydefault}{\mddefault}{\updefault}{\color[rgb]{0,0,0}$-6$}%
}}}}
\put(4096,-5596){\makebox(0,0)[lb]{\smash{{\SetFigFontNFSS{9}{10.8}{\familydefault}{\mddefault}{\updefault}{\color[rgb]{0,0,0}$-5$}%
}}}}
\put(4231,-6721){\makebox(0,0)[lb]{\smash{{\SetFigFontNFSS{9}{10.8}{\familydefault}{\mddefault}{\updefault}{\color[rgb]{0,0,0}$-2$}%
}}}}
\put(4276,-7081){\makebox(0,0)[lb]{\smash{{\SetFigFontNFSS{9}{10.8}{\familydefault}{\mddefault}{\updefault}{\color[rgb]{0,0,0}$-1$}%
}}}}
\put(3151,-2806){\makebox(0,0)[lb]{\smash{{\SetFigFontNFSS{9}{10.8}{\familydefault}{\mddefault}{\updefault}{\color[rgb]{0,0,0}$0$}%
}}}}
\put(4591,-6226){\makebox(0,0)[lb]{\smash{{\SetFigFontNFSS{9}{10.8}{\familydefault}{\mddefault}{\updefault}{\color[rgb]{0,0,0}$3$}%
}}}}
\put(5716,-4336){\makebox(0,0)[lb]{\smash{{\SetFigFontNFSS{9}{10.8}{\familydefault}{\mddefault}{\updefault}{\color[rgb]{0,0,0}$9$}%
}}}}
\put(5131,-4471){\makebox(0,0)[lb]{\smash{{\SetFigFontNFSS{9}{10.8}{\familydefault}{\mddefault}{\updefault}{\color[rgb]{0,0,0}$11$}%
}}}}
\put(9541,-3166){\makebox(0,0)[lb]{\smash{{\SetFigFontNFSS{9}{10.8}{\familydefault}{\mddefault}{\updefault}{\color[rgb]{0,0,0}$c$}%
}}}}
\put(7651,-2356){\makebox(0,0)[lb]{\smash{{\SetFigFontNFSS{9}{10.8}{\familydefault}{\mddefault}{\updefault}{\color[rgb]{0,0,0}$v_0$}%
}}}}
\put(10306,-3301){\makebox(0,0)[lb]{\smash{{\SetFigFontNFSS{9}{10.8}{\familydefault}{\mddefault}{\updefault}{\color[rgb]{0,0,0}$F_{\infty}$}%
}}}}
\put(4996,-5461){\makebox(0,0)[lb]{\smash{{\SetFigFontNFSS{9}{10.8}{\familydefault}{\mddefault}{\updefault}{\color[rgb]{0,0,0}$1$}%
}}}}
\put(10081,-5551){\makebox(0,0)[lb]{\smash{{\SetFigFontNFSS{9}{10.8}{\rmdefault}{\mddefault}{\updefault}{\color[rgb]{0,0,0}$\mathbf{1}$}%
}}}}
\put(9406,-6451){\makebox(0,0)[lb]{\smash{{\SetFigFontNFSS{9}{10.8}{\rmdefault}{\mddefault}{\updefault}{\color[rgb]{0,0,0}$\mathbf{2}$}%
}}}}
\put(9406,-7306){\makebox(0,0)[lb]{\smash{{\SetFigFontNFSS{9}{10.8}{\rmdefault}{\mddefault}{\updefault}{\color[rgb]{0,0,0}$\mathbf{1}$}%
}}}}
\end{picture}%
\end{center}
\caption{\'Etape 2: num\'erotation de quelques coins dans $F_{\infty}$
  et construction d'une corde entre les deux c\^ot\'es de la colonne.}
\label{Step2}
\end{figure}

La deuxi\`eme \'etape de la construction se d\'eroule de mani\`ere
ind\'ependante dans chacune des faces de $\mathcal{M}_0$, $F_\infty$
comprise. Soit $F$ une face de $\mathcal{M}_0$ et $c_0$ son coin
d'\'etiquette $0$. Si $F$ contient un nombre fini de sommets, donc un
nombre fini de coins, on num\'erote ses coins de $0$ \`a
$k-1$ en parcourant la face dans le sens des aiguilles d'une montre,
et en commen\c{c}ant par $c_0$. Si $F$ est la face infinie, les coins \`a droite de
la colonne sont num\'erot\'es par les entiers positifs en parcourant la
fronti\`ere dans le sens des aiguilles d'une montre et en commen\c{c}ant
par $c_0$ de num\'ero $0$. Les coins \`a gauche de la colonne sont 
num\'erot\'es par les entiers n\'egatifs en parcourant la
fronti\`ere dans le sens inverse des aiguilles d'une montre et en commen\c{c}ant
par $c_0$ de num\'ero $0$ (voir figure \ref{Step2}). On note aussi $\ell(i)$ l'\'etiquette du coin
num\'erot\'e $i$, de telle sorte que $\ell(0) = 0$ et $\ell(1) = \ell(k -
1) = 1$ dans le cas fini; $\ell(1) = \ell (-1) = 1$ dans le cas infini
(la fonction $\ell$ d\'epend \'evidemment de la face consid\'er\'ee).

Pour chaque face, finie ou infinie, on d\'efinit maintenant une fonction
successeur $s$ pour tous les coins except\'e les coins d'\'etiquettes
$0$ ou $1$ par:
\[ s(i) = \begin{cases}
\min \left\{ j > i : \, \ell(j) = \ell(i) - 1 \right\} & \text{si $i < 0$},\\
\min \left\{ j > i : \, \ell(j) = \ell(i) - 1 \right\} & \text{si $i>0$ et
$ \left\{ j>i : \, \ell(j) = \ell(i) -1 \right\} \neq \emptyset $},\\
\min \left\{ j \leqslant 0 : \, \ell(j) = \ell(i) - 1 \right\} & \text{si
  $i> 0$ et $ \left\{ j>i : \, \ell(j) = \ell(i) -1 \right\} = \emptyset $.}
\end{cases} \]
Dans le cas o\`u la face est finie, seul le deuxi\`eme cas se
produit. Dans le cas de la face infinie, la derni\`ere propri\'et\'e
de la d\'efinition \ref{colonneens} est utilis\'ee pour dire que $\{ j
\leqslant 0 : \, \ell(j) = \ell(i) - 1\}$ est fini.

\bigskip

\begin{etape} Dans chaque face, pour chaque
coin $i$ d'\'etiquette $\ell(i) \geqslant 2$, si
$|s(i) - i| \neq 1$, on ajoute une corde $(i, s(i))$ \`a l'int\'erieur de la face.
\end{etape}

\begin{proposition}[\cite{CD}, Property 6.1]
\label{intersect}
L'\'etape 2 peut se faire de mani\`ere \`a ce que les diff\'erentes cordes
$(i, s(i))$ ne se croisent pas.
\end{proposition}

\begin{remarque}
La condition $|s(i) - i| \neq 1$ revient \`a dire que la corde $(i,s(i))$
n'existe pas d\'ej\`a dans $\omega$. Dans le cas
de $F_{\infty}$, il peut arriver qu'une corde $(i, s(i))$
relie deux sommets qui sont de part et d'autre de la colonne (voir la
figure \ref{Step2}, cela se produit dans le troisi\`eme cas de la
d\'efinition de $s(i)$, ce qui correspond au cas o\`u, dans le
contour du c\^ot\'e droit de la colonne, le coin $i$ appara\^it apr\`es la
derni\`ere occurrence de l'\'etiquette $\ell(i) - 1$).
\end{remarque}

\bigskip

Cette 2\ieme\ \'etape d\'efinit une carte planaire r\'eguli\`ere
$\mathcal{M}_1$ unique \`a hom\'eomorphisme conservant l'orientation pr\`es,
et dont les faces sont d\'ecrites dans la proposition suivante:

\begin{proposition}[\cite{CD}, Property 6.2]
\label{faces}
Les faces de $\mathcal{M}_1$ sont soit triangulaires et de sommets
d'\'etiquettes $l$, $l+1$, $l+1$; soit quadrangulaires et de sommets
d'\'etiquettes $l$, $l+1$, $l+2$, $l+1$.
\end{proposition}

\bigskip

\begin{etape} On supprime toutes les ar\^etes de $\mathcal{M}_1$
reliant deux sommets de m\^eme \'etiquette.
\end{etape}
Alors (voir \cite{CD} pour les détails), on obtient une
quadrangulation infinie régulière, notée $\bij(\omega)$. De plus les
étiquettes des sommets de $\omega$ sont les distances à la racine des
sommets correspondants dans $\bij(\omega)$. Enfin, l'application
$\bij$ ainsi obtenue est une bijection de $\spinelab$ sur $\bij (\spinelab)$.

\section{Quadrangulations infinies uniformes}
\label{sec:quadrangulations}

Nous allons voir dans cette section deux mani\`eres possibles de d\'efinir
une quadrangulation infinie al\'eatoire de loi uniforme, qui donnent a
priori deux objets diff\'erents.

\subsection{D\'efinition directe}

Cette m\'ethode, introduite dans \cite{Kr}, consiste \`a d\'efinir la
quadrangulation infinie uniforme comme la loi limite de quadrangulations
al\'eatoires finies. Pour cela nous devons dans un premier temps munir
$\carte[]$, l'ensemble des quadrangulations enracin\'ees finies, d'une topologie.
Si $q \in \carte[]$, on note
$\boule{\carte[]}(q)$ l'union des faces de $q$ dont au moins un des sommets est
\`a distance strictement inf\'erieure \`a $R$ de la racine
($\boule{\carte[]}(q)$ est donc une carte planaire enracin\'ee
finie). On munit $\carte[]$ de la distance :
\[ \dist{\carte[]} (q_1,q_2) = \left( 1 + \sup \left\{ R : \,
\boule{\carte[]}(q_1) = \boule{\carte[]}(q_2) \right\} \right)^{-1} ,\]
l'\'egalit\'e \'etant au sens de l'\'egalit\'e entre deux cartes planaires enracin\'ees finies.

Soit $(\kcartes,\dist{\kcartes})$ le compl\'et\'e de
$(\carte[],\dist{\carte[]})$. Les \'el\'ements de $\kcartes$ 
autres que les quadrangulations finies sont appel\'es quadrangulations enracin\'ees
infinies au sens de Krikun. Cette d\'efinition n'est pas \'equivalente \`a celle donn\'ee
dans la d\'efinition \ref{Qinf}.

\begin{theoreme}[\cite{Kr}, Theorem 1]
Pour tout $n \geqslant 1$ soit $\mcartes$ la mesure de
probabilit\'e uniforme sur $\carte[n]$. La suite $(\mcartes)_{n \in \mathbb{N}}$
converge faiblement dans l'espace $(\kcartes, \dist{\kcartes})$ vers une mesure
de probabilit\'e not\'ee $\mcartes[]$. De plus, $\mcartes[]$ est port\'ee par
l'ensemble des quadrangulations enracin\'ees infinies (au sens de Krikun).
\end{theoreme}

\begin{remarque}
On peut prolonger la fonction $q \in \carte[] \mapsto \boule{\carte[]}(q)$
en une fonction continue $\boule{\kcartes}$ sur $\kcartes$.
$\boule{\kcartes} (q)$ s'interpr\`ete de mani\`ere naturelle comme
l'union des faces de $q$ dont au
moins un des sommets est \`a distance strictement inf\'erieure \`a $R$ de
la racine.
\end{remarque}

\subsection{D\'efinition indirecte}
\label{sec:indirect}

Une autre approche possible pour d\'efinir la quadrangulation infinie
uniforme est de d\'efinir dans un premier temps l'arbre bien \'etiquet\'e
infini uniforme et de consid\'erer la mesure image de sa loi par la bijection
de Schaeffer. C'est l'approche choisie dans l'article \cite{CD} auquel nous
renvoyons pour plus de d\'etails et les preuves des r\'esultats qui
vont suivre. On munit l'ensemble $\arbres$ des arbres bien étiquetés
(finis ou infinis) de la distance suivante:
\[ \dist{\arbres}(\omega,\omega') = \left( 1 + \sup \left\{ S : \,
  \boule[S]{\arbres}(\omega) = \boule[S]{\arbres} (\omega') \right\} \right)^{-1} ,\]
o\`u $\boule[S]{\arbres} (\omega)$ est le sous arbre de $\omega$
jusqu'\`a la g\'en\'eration $S$. L'espace métrique $\left(\arbres, \dist{\arbres}\right)$
est complet.

On a alors le r\'esultat suivant:
\begin{theoreme}[\cite{CD}, Theorem 3.1]
\label{cvarbres}
La suite $(\marbres)_{n \in \mathbb{N}}$ des probabilit\'es
uniformes sur les arbres bien \'etiquet\'es \`a $n$ ar\^etes converge faiblement vers
une probabilit\'e $\marbres[]$ \`a support dans
$\arbre$. On appelle cette loi limite la loi de l'arbre bien
\'etiquet\'e infini uniforme.
\end{theoreme}

Une des \'etapes cl\'e de la preuve de ce r\'esultat est la preuve de la convergence
\[ \marbres \left( \omega \in \arbres : \boule[S]{\arbres}
(\omega) = \omega^{\star} \right) \underset{n \to \infty}{\longrightarrow} \marbres[] \left(
\omega \in \arbres : \boule[S]{\arbres} 
(\omega) = \omega^{\star} \right) \]
pour $S >0$ et $\omega^{\star}$ un arbre bien \'etiquet\'e fini de hauteur
$S$. Ceci est fait au moyen de calculs explicites:
pour tout $\omega^{\star}$ bien \'etiquet\'e de hauteur $S$
ayant exactement $k$ sommets de g\'en\'eration $S$ d'\'etiquettes
respectives $l_1, \ldots , l_k$:
\begin{align}
\label{boulen}
\marbres \left( \omega \in \arbres : \boule[S]{\arbres}
(\omega) = \omega^{\star} \right) & = \frac{1}{D_n} \sum_{n_1 + \cdots + n_k
= n - |\omega ^{\star}|} \prod_{j=1}^k D_{n_j}^{(l_j)} ,\\
\label{boule}
\marbres[] \left( \omega \in \arbres : \boule[S]{\arbres}
(\omega) = \omega^{\star} \right) & = \frac{1}{12^{|\omega^{\star}|}} \sum_{i =
1}^k d_{l_i} \prod_{j \neq i} w_{l_j},
\end{align}
o\`u $D_n^{(l)}$ est le cardinal de $\arbre[n]^{(l)}$ pour tout
$l \geqslant 1$,
$D_n^{(1)} = D_n$ et
\begin{align}
\label{defw}
w_{l} & = 2 \frac{l (l +3)}{(l +1)(l +2)},\\
d_{l} & = \frac{2w_{l}}{560} (4 l^4 + 30 l^3 + 59 l^2 + 42 l
+ 4),
\label{defd}
\end{align}
pout tout $l \geqslant 1$.

\begin{proposition}[\cite{CD}, Theorem 5.9]
\label{etiquettes}
$\marbres[]$ est \`a support dans $\spine$,
  l'ensemble des arbres bien \'etiquet\'es infinis ayant une unique
  colonne. De plus $\marbres[]$ a la propri\'et\'e
  suivante :
\[\mathbb{E}_{\marbres[]} \left[ N_l(\omega) \right] =
  {\it O}(l^3) \quad \text{quand $l \to \infty$}.\]
Ainsi, la loi de l'arbre bien \'etiquet\'e infini uniforme est \`a support
dans $\spinelab$, l'ensemble des arbres bien \'etiquet\'es infinis ayant
une seule colonne et dont toutes les \'etiquettes n'apparaissent qu'un
nombre fini de fois.
\end{proposition}

Un arbre de loi $\marbres[]$ a donc p.s. une unique colonne. \cite{CD}
donne une description pr\'ecise du processus des \'etiquettes des
sommets de cette colonne et des sous-arbres attach\'es \`a chacun de
ces sommets. Cette description fait intervenir la mesure $\rho^{(l)}$
d\'efinie sur $\mathbf{T}^{(l)}$ par $\rho^{(l)} (\omega) =
12^{-|\omega|}$ pour tout $\omega \in \mathbf{T}^{(l)}$. On v\'erifie
facilement que $\frac{1}{2} \rho^{(l)}$ est la loi de l'arbre de
Galton-Watson de loi de reproduction g\'eom\'etrique de param\`etre
$\frac{1}{2}$, muni d'\'etiquettes al\'eatoires selon les r\`egles
suivantes: la racine a pour \'etiquette $l$ et chaque sommet autre que
la racine a une \'etiquette choisie uniform\'ement dans $\{m-1,m,m+1
\}$ o\`u $m$ est l'\'etiquette de son p\`ere, ces choix \'etant faits
ind\'ependamment pour tous les sommets. La proposition 2.4 de
\cite{CD} montre que $\rho^{(l)}
(\arbre[]^{(l)}) = w_l$. Cela permet d'introduire la mesure de probabilit\'e
$\hat{\rho}^{(l)}$ sur $\arbre[]^{(l)}$ d\'efinie par
$\hat{\rho}^{(l)} (\omega) = w_l^{-1} \rho^{(l)} (\omega) = w_l^{-1}
12^{-|\omega|}$ pour tout $\omega \in \arbre[]^{(l)}$. On a le
th\'eor\`eme suivant:
\begin{theoreme}[\cite{CD}, Theorem 4.4]
\label{descmu}
Soit $\omega$ un arbre bien \'etiquet\'e de loi $\mu$, et soit $u_0,
u_1, u_2, \ldots $ la suite des sommets de la colonne de $\omega$,
énumérés dans l'ordre généalogique. Pour tout $n \geqslant 0$, soit
$Y_n$ l'étiquette de $u_n$.
\begin{enumerate}
\item Le processus $(Y_n)_{n \geqslant 0}$ est une cha\^ine de Markov \`a
valeurs dans $\mathbb{N}$ de noyau de transition $\Pi$ d\'efini par:
\begin{align*}
\Pi (l,l-1) & = q_l = \frac{(w_l)^2}{12 d_l} d_{l-1} & \text{si $l
\geqslant 2$,}\\
\Pi (l,l) & = r_l = \frac{(w_l)^2}{12} & \text{si $l \geqslant 1$,}\\
\Pi (l,l+1) & = p_l = \frac{(w_l)^2}{12 d_l} d_{l+1} & \text{si $l \geqslant 1$.}
\end{align*}
\item Conditionnellement \`a $(Y_n)_{n \geqslant 0} = (y_n)_{n
\geqslant 0}$, la suite $(L_n)_{n \geqslant 0}$ des sous-arbres de $\omega$
attach\'es \`a gauche de chaque sommet de la colonne et la suite
$(R_n)_{n \geqslant 0}$ des sous-arbres attach\'es \`a droite sont
ind\'ependantes et constitu\'ees d'arbres ind\'ependants entre eux de
lois respectives $\hat{\rho}^{(y_n)}$.
\end{enumerate}
\end{theoreme}

Nous allons maintenant transporter la loi de l'arbre bien \'etiquet\'e
infini uniforme sur les quadrangulations au moyen de la bijection de
Schaeffer. Munissons $\bij ( \spinelab ) $ de la distance $\dist{\bij}$ de telle sorte
que $\bij$ soit une isom\'etrie entre $\spinelab$ et $\bij(\spinelab)$.
On note $\mimage$ et $\mimage[]$ les mesures images de
$\marbres$ et $\marbres[]$ par
$\bij$. La mesure $\mimage[]$ est bien d\'efinie car
$\marbres[]$ est \`a support dans $\spinelab$.

Comme $\bij$ met en bijection $\arbre[n]$ et $\carte[n]$,
$\mimage = \mcartes$ est la probabilit\'e uniforme sur les
quadrangulations \`a $n$ faces. Une cons\'equence
directe du th\'eor\`eme \ref{cvarbres} est alors que,
dans $\left( \bij(\spinelab), \dist{\bij} \right)$, $(\mimage)_{N \in \mathbb{N}}$
converge faiblement vers $\mimage[]$,
peut aussi \^etre consid\'er\'ee comme
une probabilit\'e uniforme sur l'ensemble des quadrangulations infinies.

\begin{remarque}
La topologie associ\'ee \`a $\dist{\bij}$ n'est pas la m\^eme que
celle donn\'ee pr\'ec\'edemment pour $\kcartes$, il est donc naturel
de se demander si les deux notions de quadrangulation infinie uniforme
que nous avons donn\'ees co\"incident.
\end{remarque}

\section{\'Egalit\'e des deux quadrangulations infinies uniformes}
\label{sec:egalite}

Notre but est maintenant de montrer que les deux notions de quadrangulation infinie
uniforme que nous avons d\'efinies sont les m\^emes. La premi\`ere question
que l'on doit se poser concerne les liens entre les deux d\'efinitions  des
quadrangulations infinies, \`a savoir
$\bij(\spinelab)$ qui est un sous-ensemble de l'ensemble des
quadrangulations consid\'er\'ees comme cartes planaires r\'eguli\`eres
d'une part, et $\kcartes$ le compl\'et\'e de $\carte[]$ d'autre part. Nous allons voir que
$\bij(\spinelab)$ s'identifie \`a un sous-ensemble de $\kcartes$, ce qui
permet de consid\'erer $\mimage[]$ comme une mesure sur
$\kcartes$, port\'ee par $\bij(\spinelab)$.

Plus pr\'ecis\'ement, soient $R>0$ et $\omega \in \spinelab$. On d\'efinit
$B_R(\bij(\omega))$ comme l'union des faces de $\bij(\omega)$ contenant au
moins un sommet \`a distance strictement plus petite que $R$ de
l'origine. Comme $\omega$ n'a qu'un nombre fini de sommets d'\'etiquette
inf\'erieure ou \'egale \`a $R+1$, il n'y a qu'un nombre fini de telles faces et
$B_R(\bij(\omega))$ est une carte finie. Cette carte est donc telle
que $\mathbb{S}^2 \setminus B_R(\bij(\omega))$ a un
nombre fini de composantes connexes, et les fronti\`eres de ces composantes
sont des cycles de $\bij(\omega)$ de longueur finie.

Soit $\gamma$ un tel cycle. Chaque
ar\^ete constituant $\gamma$ est commune \`a une face de $\bij(\omega)$
ayant au moins un sommet \`a distance strictement plus petite que $R$ de
l'origine et \`a une autre face n'ayant que des sommets \`a distance au moins
$R$ de l'origine. Le fait que la quadrangulation soit bipartie nous assure alors que chaque
ar\^ete de $\gamma$ est une ar\^ete $(R,R+1)$. Ainsi en ajoutant
\`a $B_R(\bij(\omega))$ un sommet \`a l'int\'erieur
de $\gamma$ que l'on relie aux sommets de $\gamma$ \`a distance $R+1$ de la
racine, et cela pour chaque cycle $\gamma$ fronti\`ere d'une composante connexe de 
$\mathcal{S}^2 \setminus B_R(\bij(\omega))$, on obtient
une quadrangulation finie. La suite des quadrangulations finies ainsi
construites pour chaque $R > 0$ converge au sens de la
topologie de Krikun vers $\bij(\omega)$ quand $R$ tend vers l'infini. Cela montre que pour
tout $\omega \in \spinelab$, $\bij(\omega)$ s'identifie \`a un
\'el\'ement de $\kcartes$.

Pour pouvoir consid\'erer $\mimage[]$ comme une mesure sur $\kcartes$, il
suffit maintenant de voir que $\bij : \spinelab \rightarrow \kcartes$ est
mesurable pour la tribu bor\'elienne associ\'ee \`a
$\dist{\kcartes}$. On a pour cela besoin du lemme suivant, prouv\'e
dans la section \ref{propbij}:
\begin{lemme}
\label{mesurable}
Soit $R > 0$ et $\omega_0 \in \spinelab$. L'ensemble $ A = \{ \omega
\in \spinelab : \, \boule{\kcartes} (\bij(\omega)) =
\boule{\kcartes} (\bij(\omega_0)) \}$ est mesurable pour la
tribu bor\'elienne associ\'ee \`a $\dist{\arbres}$.
\end{lemme}

Soit $q^{\star} \in \kcartes$, alors d'apr\`es ce dernier lemme
\begin{equation*}
 \bij^{-1} \left( \left\{ q \in
    \kcartes : \, \dist{\kcartes} (q,q^{\star}) \leqslant \frac{1}{R+1} \right\} \right)
= \bij^{-1} \left( \left\{ \vphantom{\frac{1}{R}} q \in \kcartes : \, \boule{\kcartes} (q) =
    \boule{\kcartes} (q^{\star}) \right\} \right)
\end{equation*}
est mesurable pour la tribu bor\'elienne associ\'ee \`a
  $\dist{\arbres}$. Cela suffit pour obtenir que
$\bij : \spinelab \rightarrow \kcartes $ est mesurable.

Nous pouvons maintenant énoncer notre résultat principal.

\begin{theoreme}
\label{egalite}
La suite $(\mimage)_{n \in \mathbb{N}}$ converge faiblement vers
$\mimage[]$ au sens de Krikun, i.e. dans $\left(
\kcartes , \dist{\kcartes} \right)$. En conséquence, $\mimage[]$, vue
comme mesure de probabilité sur $\left(\kcartes , \dist{\kcartes}
\right)$, coïncide avec $\mcartes[]$.
\end{theoreme}
La seconde assertion est une conséquence immédiate de la première
puisque $\mimage = \mcartes$ et $\mcartes[]$ est définie comme la
limite au sens de Krikun de la suite $(\mcartes)$.

Pour établir la première assertion, il suffira de montrer que pour
tous $q^{\star} \in \kcartes$ et $R> 0$ on a
\begin{equation*}
\mimage \Big( q \in \kcartes : 
\boule{\kcartes}(q) = \boule{\kcartes}(q^{\star})
\Big) \underset{n \to \infty}{\longrightarrow}
\mimage[] \Big( q \in \kcartes :
\boule{\kcartes}(q) = \boule{\kcartes}(q^{\star}) \Big)
\end{equation*}
o\`u on a identifi\'e $\bij(\spinelab)$ \`a un sous-ensemble de
$\kcartes$.

Le reste de ce travail est consacré à la preuve de cette convergence.

\subsection{Une propri\'et\'e de la correspondance de Schaeffer}
\label{propbij}

Pour tous $S>0$ et $R>0$ on note
\begin{multline}
\label{omegas}
\Omega_S (R) = \{ \omega \in \arbres : \,
\omega \text{ a un sommet d'\'etiquette} \leqslant R+1 \\ \text{
  strictement au dessus de la génération } S \}.
\end{multline}
Dans les deux premiers r\'esultats de cette partie, $S$ et $R$ sont
deux entiers strictement positifs fix\'es.

\begin{proposition}
\label{facegen}
Si $\omega$ est un arbre de $\spinelab$ n'appartenant pas \`a
$\Omega_S(R)$ (i.e. tel que tout sommet de g\'en\'eration strictement
plus grande que $S$ est d'\'etiquette au moins \'egale \`a $R+2$), alors
$\boule{\kcartes} (\bij (\omega)) =
\boule{\kcartes} (\bij (\boule[S]{\arbres}(\omega)))$.
\end{proposition}
\begin{proof}
Nous allons reprendre la construction de $\bij(\omega)$ \'etape par
\'etape comme dans la section \ref{sch}. On fixe un plongement de
$\omega$ dans $\mathbb{S}^2$ comme carte planaire accepatable.

Pendant la premi\`ere \'etape, on construit \`a partir de $\omega$ une nouvelle
carte planaire acceptable $\mathcal{M}_0(\omega)$ en reliant un sommet $v_0$
d'\'etiquette $0$ aux coins d'\'etiquette $1$. Construisons de m\^eme
$\mathcal{M}_0( \boule[S]{\arbres} (\omega))$. Les ar\^etes ajout\'ees sont
uniquement d\'etermin\'ees par les
coins d'\'etiquette $1$, coins donn\'es par la structure de
$\boule[S]{\arbres} (\omega)$ (en effet aucun sommet au dessus de la
g\'en\'eration $S$ n'est d'\'etiquette inf\'erieure ou \'egale \`a $2$). Les deux
cartes dont nous disposons ont donc le m\^eme nombre $p$ de faces, qui ont
de plus les m\^emes fronti\`eres constitu\'ees chacune par les deux ar\^etes joignant
$v_0$ aux deux coins d'\'etiquette $1$ ainsi que, dans le cas o\`u ces coins
d'\'etiquette $1$ appartiennent \`a deux sommets distincts, les ar\^etes
figurant sur l'unique chemin injectif joignant ces deux
sommets d'\'etiquette $1$ dans l'arbre.
Notons $F_1(\omega) , \ldots ,
F_p(\omega)$ et $F_1 (\boule[S]{\arbres}(\omega)), \ldots ,
F_p(\boule[S]{\arbres}(\omega))$ ces faces, de mani\`ere \`a ce que
pour tout $i$ les fronti\`eres de $F_i(\omega)$ et $F_i
(\boule[S]{\arbres}(\omega))$ soient identiques.

\bigskip

La deuxi\`eme \'etape consiste \`a construire pour chaque coin $c$ une ar\^ete $(c, s(c))$
de mani\`ere ind\'ependante dans chaque face afin d'obtenir deux cartes
planaires r\'eguli\`eres que nous noterons $\mathcal{M}_1(\omega)$ et
$\mathcal{M}_1(\boule[S]{\arbres}(\omega))$. Consid\'erons donc $F_i(\omega)$ une
face de $\mathcal{M}_0(\omega)$ et $F_i(\boule[S]{\arbres}(\omega))$ la face de
$\mathcal{M}_0( \boule[S]{\arbres} (\omega))$ correspondante. 
Les coins de ces faces sont num\'erot\'es $(c_{i,j})_{j \in J_i}$ pour
$F_i(\omega)$ et $(c'_{i,j})_{j \in J'_i}$ pour
$F_i(\boule[S]{\arbres}(\omega))$, la num\'erotation correspondant \`a
l'ordre dans lequel on rencontre les coins en parcourant les faces dans
le sens des aiguilles d'une montre si elles sont finies, et en
parcourant la partie \`a gauche de la colonne dans le sens inverse des
aiguilles d'une montre, la partie \`a droite de la colonne dans le
sens des aiguilles d'une montre dans le cas de la face infinie.

Pour chaque $i \in \{1, \ldots , p \}$,
soient $v_{i,1}, \ldots , v_{i,k_i}$ les sommets de g\'en\'eration $S$ de
$F_i(\omega)$ ayant au moins un descendant (ce sont aussi des sommets
de g\'en\'eration $S$ de $F_i(\boule[S]{\arbres}(\omega))$), ces
sommets sont d'\'etiquette au moins \'egale \`a $R+1$.
Pour tout $j \leqslant k_i$, soit $e_{i,j}$ le
dernier coin d'\'etiquette $R+1$ avant $v_{i,j}$ dans
$F_i(\omega)$. Ce coin est identique dans
$F_i(\omega)$ et $F_i(\boule[S]{\arbres}(\omega))$. On construit donc dans
$F_i(\omega)$ et $F_i(\boule[S]{\arbres}(\omega))$ la m\^eme ar\^ete
$(e_{i,j}, s(e_{i,j}))$ reliant $e_{i,j}$ au premier coin
d'\'etiquette $R$ suivant $e_{i,j}$ (ce coin est aussi le premier coin
d'\'etiquette $R$ suivant tous les coins de $v_{i,j}$, voir la figure \ref{faceisole}).

\begin{figure}[!t]
\begin{center}
\begin{picture}(0,0)%
\includegraphics{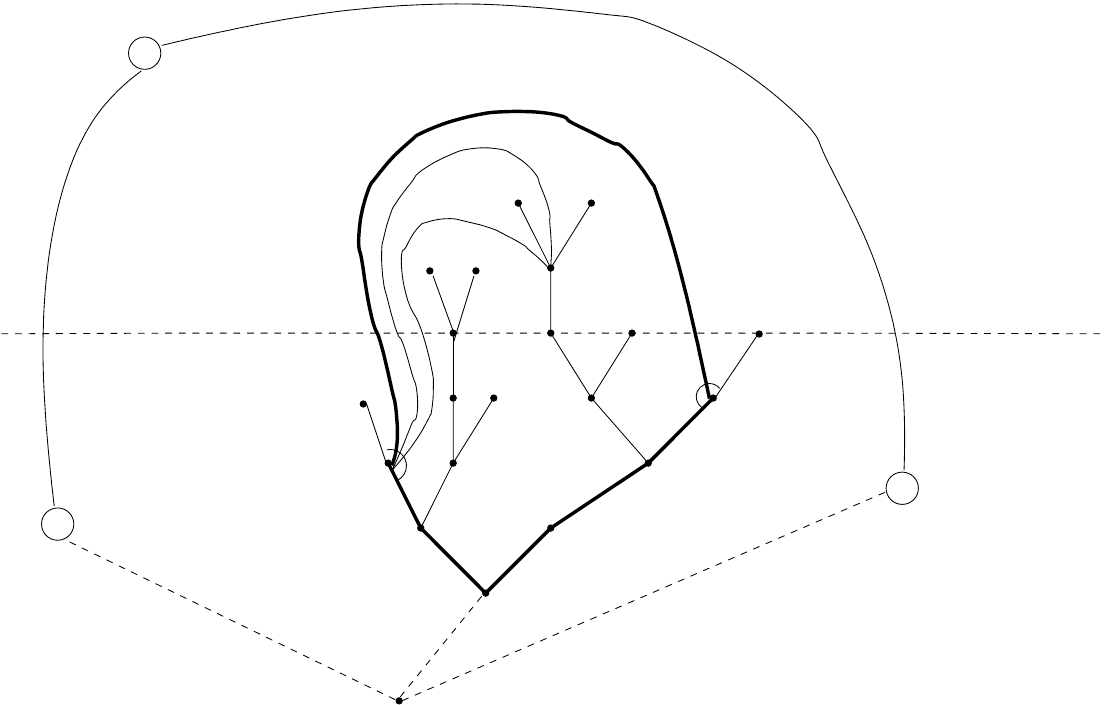}%
\end{picture}%
\setlength{\unitlength}{2155sp}%
\begingroup\makeatletter\ifx\SetFigFontNFSS\undefined%
\gdef\SetFigFontNFSS#1#2#3#4#5{%
  \reset@font\fontsize{#1}{#2pt}%
  \fontfamily{#3}\fontseries{#4}\fontshape{#5}%
  \selectfont}%
\fi\endgroup%
\begin{picture}(9807,6174)(874,-5858)
\put(5401,-2806){\makebox(0,0)[lb]{\smash{{\SetFigFontNFSS{6}{7.2}{\familydefault}{\mddefault}{\updefault}{\color[rgb]{0,0,0}$v_{i,j}$}%
}}}}
\put(7246,-3301){\makebox(0,0)[lb]{\smash{{\SetFigFontNFSS{6}{7.2}{\familydefault}{\mddefault}{\updefault}{\color[rgb]{0,0,0}$e_{i,j}$}%
}}}}
\put(6841,-1411){\makebox(0,0)[lb]{\smash{{\SetFigFontNFSS{6}{7.2}{\rmdefault}{\mddefault}{\updefault}{\color[rgb]{0,0,0}$\gamma_{i,j}$}%
}}}}
\put(3511,-3931){\makebox(0,0)[lb]{\smash{{\SetFigFontNFSS{6}{7.2}{\familydefault}{\mddefault}{\updefault}{\color[rgb]{0,0,0}$s(e_{i,j})$}%
}}}}
\put(5401,-3706){\makebox(0,0)[lb]{\smash{{\SetFigFontNFSS{6}{7.2}{\familydefault}{\mddefault}{\updefault}{\color[rgb]{0,0,0}$F_{i,j}$}%
}}}}
\put(2296,-1366){\makebox(0,0)[lb]{\smash{{\SetFigFontNFSS{6}{7.2}{\rmdefault}{\mddefault}{\updefault}{\color[rgb]{0,0,0}$F_i$}%
}}}}
\put(2071,-196){\makebox(0,0)[lb]{\smash{{\SetFigFontNFSS{6}{7.2}{\rmdefault}{\mddefault}{\updefault}{\color[rgb]{0,0,0}$\mathbf{0}$}%
}}}}
\put(10666,-2581){\makebox(0,0)[lb]{\smash{{\SetFigFontNFSS{6}{7.2}{\familydefault}{\mddefault}{\updefault}{\color[rgb]{0,0,0}g\'en\'eration $S$}%
}}}}
\put(5806,-2131){\makebox(0,0)[lb]{\smash{{\SetFigFontNFSS{6}{7.2}{\rmdefault}{\mddefault}{\updefault}{\color[rgb]{0,0,0}$u$}%
}}}}
\put(4906,-1366){\makebox(0,0)[lb]{\smash{{\SetFigFontNFSS{6}{7.2}{\rmdefault}{\mddefault}{\updefault}{\color[rgb]{0,0,0}$f$}%
}}}}
\put(1351,-4336){\makebox(0,0)[lb]{\smash{{\SetFigFontNFSS{6}{7.2}{\rmdefault}{\mddefault}{\updefault}{\color[rgb]{0,0,0}$\mathbf{1}$}%
}}}}
\put(8776,-4021){\makebox(0,0)[lb]{\smash{{\SetFigFontNFSS{6}{7.2}{\rmdefault}{\mddefault}{\updefault}{\color[rgb]{0,0,0}$\mathbf{1}$}%
}}}}
\end{picture}%
\end{center}
\caption{Une face $F_i$ et un cycle $\gamma_{i,j}$ construit \`a
  partir d'un sommet $v_{i,j}$ de g\'en\'eration $S$.}
\label{faceisole}
\end{figure}
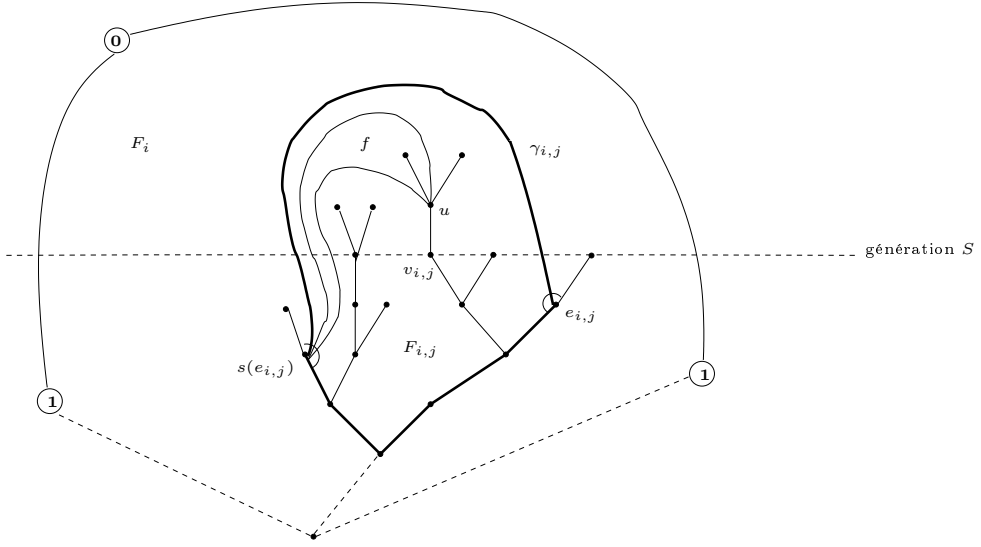

Dans $F_i(\omega)$ et dans $F_i(\boule[S]{\arbres}(\omega))$, on
obtient ainsi, pour chaque $j \in \{ 1, \ldots , k_i\}$, le m\^eme
cycle $\gamma_{i,j}$ constitué de l'ar\^ete $(e_{i,j}, s(e_{i,j}))$ et
du chemin g\'en\'ealogique entre $e_{i,j}$ et $s(e_{i,j})$ (voir la
figure \ref{faceisole}). Si $j \neq j'$, ou bien les int\'erieurs des
cycles $\gamma_{i,j}$ et $\gamma_{i,j'}$ sont disjoints, ou bien
l'int\'erieur de $\gamma_{i,j}$ est contenu dans l'int\'erieur de
$\gamma_{i,j'}$ ou vice-versa. Ici, l'int\'erieur (strict) d'un
cycle est d\'efini comme la composante connexe du compl\'ementaire du
cycle qui ne contient pas $v_0$.

Montrons maintenant que toute face $f$ de $\Phi(\omega)$ qui rencontre
l'intérieur d'un cycle $\gamma_{i,j}$ ne contient que des sommets
d'\'etiquette au moins \'egale \`a $R$. Traitons d'abord le cas o\`u
il existe un sommet $u$ de $f$ qui appartient \`a l'int\'erieur de
$\gamma_{i,j}$. Si l'\'etiquette de $u$ est sup\'erieure ou \'egale \`a
$R+2$, le r\'esultat est \'evident. Sinon l'\'etiquette de $u$ est
$R+1$ et, pour que la face $f$ contienne un sommet d'\'etiquette
$R-1$, il serait n\'ecessaire que $u$ soit reli\'e par deux ar\^etes
\`a des sommets d'\'etiquette $R$: le seul sommet d'\'etiquette $R$
possible est $s(e_{i,j})$, et le dernier sommet de la face $f$ devrait
\^etre \`a l'int\'erieur de la r\'eunion des deux ar\^etes joignant
$u$ \`a $s(e_{i,j})$ (cf. figure \ref{faceisole}), ce qui rend
impossible le fait qu'il soit d'\'etiquette $R-1$. Le cas o\`u $f$
rencontre l'int\'erieur du cycle sans qu'aucun de ses sommets ne soit
dans cet int\'erieur est trait\'e par des arguments analogues.

\bigskip

La discussion pr\'ec\'edente montre que les faces de $\Phi(\omega)$,
respectivement de $\Phi \left( B_{\mathbb{T},S} (\omega) \right)$, qui
rencontrent l'int\'erieur d'un cycle $\gamma_{i,j}$ ne sont pas prises
en compte dans la d\'efinition de la boule
$B_{\overline{\mathbf{Q}},R} \left( \Phi (\omega) \right)$,
respectivement $B_{\overline{\mathbf{Q}},R} \left( \Phi \left(
B_{\mathbb{T},S} (\omega) \right) \right)$.

Notons $\Phi_1(\omega)$ la carte palanire obtenue dans la deuxi\`eme
\'etape de la construction de $\Phi(\omega)$ (cette carte est not\'ee
$\mathcal{M}_1$ dans la partie \ref{sch}). Consid\'erons alors la
carte planaire $\widetilde{\Phi}_1(\omega)$ obtenue \`a partir de
$\Phi_1(\omega)$ en supprimant les sommets se trouvant \`a
l'int\'erieur d'un des cycles $\gamma_{i,j}$ ainsi que les ar\^etes
trac\'ees \`a l'int\'erieur d'un tel cycle dans la deuxi\`eme \'etape
de la construction de la partie \ref{sch}. On a
\begin{equation}
\label{phi1}
\widetilde{\Phi}_1(\omega) = \widetilde{\Phi}_1 \left(
B_{\overline{\mathbb{T}},S} (\omega) \right)
\end{equation}
simplement parce que tous les sommets de $\omega$ au-dessus
(strictement) de la g\'en\'eration $S$ sont \`a l'int\'erieur d'un des
cycles $\gamma_{i,j}$.

Finalement, notons $\Phi_2(\omega)$ la carte obtenue \`a partir de
$\widetilde{\Phi}_1(\omega)$ en supprimant les ar\^etes de
$\widetilde{\Phi}_1(\omega)$ joignant deux sommets de m\^eme
\'etiquette inf\'erieure ou \'egale \`a $R$. D'apr\`es les
observations pr\'ec\'edentes, toute face de $\Phi(\omega)$ prise en
compte dans la boule $B_{\overline{\mathbf{Q}},R} \left( \Phi (\omega)
\right)$ est aussi une face quadrangulaire de $\Phi_2(\omega)$ et,
inversement, toute face quadrangulaire de $\Phi_2(\omega)$ contenant
un sommet d'\'etiquette strictement inf\'erieure \`a $R$ est une face
de $\Phi(\omega)$ prise en compte dans la boule
$B_{\overline{\mathbf{Q}},R} \left( \Phi (\omega) \right)$.
En d'autres mots, la boule $B_{\overline{\mathbf{Q}},R} \left( \Phi
  (\omega) \right)$ est la r\'eunion des faces quadrangulaires de
$\Phi_2(\omega)$ dont au moins un sommet est d'\'etiquette strictement
inf\'erieure \`a $R$. D'apr\`es \eqref{phi1}, on a
\[
\widetilde{\Phi}_2(\omega) = \widetilde{\Phi}_2 \left(
B_{\overline{\mathbb{T}},S} (\omega) \right)
\]
et on conclut que
\[
B_{\overline{\mathbf{Q}},R} \left( \Phi (\omega) \right) =
B_{\overline{\mathbf{Q}},R} \left( \Phi \left(
B_{\overline{\mathbb{T}},S} (\omega) \right) \right).
\]
\end{proof}

Cette propri\'et\'e permet de d\'emontrer le r\'esultat qui suit:
\begin{corollaire}
\label{pbcle}
Soit $\omega_0 \in \spinelab$.
Il existe une famille $(\omega_i^{S,R})_{i \in I}$, finie ou d\'enombrable, d'arbres
de $\spinelab \cap \Omega_S(R)^c$, v\'erifiant 
\[ \boule{\kcartes}(\bij(\omega_i^{S,R}))
  = \boule{\kcartes}(\bij(\omega_0)), \, \forall i \in I \]
  telle que, pour tout $\omega \in \spinelab \cap \Omega_S(R)^c$, les deux conditions
  suivantes sont \'equivalentes:
\begin{enumerate}
\item $\boule{\kcartes}(\bij(\omega)) =
  \boule{\kcartes}(\bij(\omega_0))$;
\item $\exists i \in I$ tel que
  $\boule[S]{\arbres}(\omega) =
  \boule[S]{\arbres}(\omega_i^{S,R})$.
\end{enumerate} 
\end{corollaire}
\begin{proof}
La famille $\left(\omega_i^{S,R}\right)_{i \in I}$ des arbres finis ayant au plus $S$
g\'en\'erations et tels que pour tout $i$ $\boule{\kcartes}(\bij(\omega_i^{S,R})) =
\boule{\kcartes}(\bij(\omega_0))$ est une famille finie ou d\'enombrable
et convient. En effet, si $\omega \in \spinelab \cap \Omega_S(R)^c$ et
s'il existe $i \in I$ tel que 
$\boule[S]{\arbres}(\omega) =
\boule[S]{\arbres}(\omega_i^{S,R})$,
alors la proposition \ref{facegen} assure que
$\boule{\kcartes}(\bij(\omega)) =
\boule{\kcartes}(\bij(\omega_i^{S,R})) =
\boule{\kcartes}(\bij(\omega_0))$. Inversement si $\boule{\kcartes}(\bij(\omega)) =
  \boule{\kcartes}(\bij(\omega_0))$ et $\omega \in \spinelab \cap \Omega_S(R)^c$
  alors $\omega' = \boule[S]{\arbres}(\omega)$ v\'erifie $\boule{\kcartes}(\bij(\omega)) =
\boule{\kcartes}(\bij(\omega'))$ d'apr\`es la proposition \ref{facegen}, et
donc $\omega'$ appartient \`a la famille $(\omega_i^{S,R})_{i \in I}$.
\end{proof}

Terminons cette section par la preuve du lemme \ref{mesurable}.
\begin{proof}[D\'emonstration du lemme \ref{mesurable}]
Soit $R > 0$ fixé. Pour $S >0$ l'ensemble $\Omega_S(R)$ est ouvert dans $\arbres$. On a de plus:
\begin{equation*}
\begin{split}
A & = \bigcup_{S > 0} \left( \left\{ \omega
\in \spinelab : \, \boule{\kcartes} (\bij(\omega)) =
\boule{\kcartes} (\bij(\omega_0)) \right\} \cap \Omega_S(R)^c \right)\\
& = \bigcup_{S > 0} \bigcup_{i \in I_{S,R}} \left( \left\{ \omega \in
  \spinelab : \, \boule[S]{\arbres} (\omega) = \boule[S]{\arbres}
(\omega_i^{S,R}) \right\} \cap \Omega_S(R)^c \right)
\end{split}
\end{equation*}
o\`u, pour tout $S$, $I_{S,R}$ est un ensemble d\'enombrable et
$\left( \omega_i^{S,R} \right)_{i \in I_{S,R}}$ est la famille
d'arbres donn\'ee par le corollaire \ref{pbcle}.
On obtient ainsi que $A$ est mesurable.
\end{proof}

\subsection{Comportement asymptotique des \'etiquettes de la colonne
  de l'arbre infini uniforme}
\label{etiquettescolonne}

Rappellons que $(Y_k)_{k \geqslant 0}$, la suite des \'etiquettes des sommets de la colonne de l'arbre
al\'eatoire infini uniforme, est une chaîne de Markov de matrice de
transition $\Pi$ introduite dans le théorème \ref{descmu}. Nous
étudions dans cette partie le comportement asymptotique de cette
cha\^ine de Markov.

\begin{lemme}
\label{labspine}
La chaîne $(Y_n)_{n \geqslant 0}$ est transitoire. De plus, 
pour tout $\varepsilon > 0$ il existe $\alpha > 0$ tel que pour tout
$k$ assez grand on ait:
\[
\mathbb{P}\left[ Y_j \geqslant \alpha k, \, \forall j \geqslant 0 \, \middle|
  \, Y_0   =k \right] \geqslant 1 - \varepsilon.
\]
\end{lemme}
\begin{proof}
On a le d\'eveloppement limit\'e $\frac{q_k}{p_k} = 1 - \frac{8}{k} +
{\it O} \left( \frac{1}{k^2} \right)$ (\cite{CD}, lemme 5.5) et donc il
existe une constante $C >0$ telle que
\[
\prod_{i=2}^{k} \frac{q_i}{p_i} \underset{k \to \infty}{\sim} C k^{-8}.
\]
Un argument standard pour les chaînes de naissance et de mort assure alors que la cha\^ine $Y$ est
transitoire et donne, pour tous $k > j \geqslant 1$,
\[
\mathbb{P}_k \left[T_j = \infty \right] =
  \frac{\sum_{i=j}^{k-1}\frac{q_i}{p_i} \cdots \frac{q_{x+1}}{p_{x+1}}}
       {\sum_{i=j}^{\infty}\frac{q_i}{p_i} \cdots \frac{q_{x+1}}{p_{x+1}}}.
\]
O\`u $T_j$ est le temps d'atteinte de $j$. On a donc, pour $\alpha <
1$ et pour tout $k$ assez grand:
\begin{align*}
\mathbb{P}_k \left[T_{[\alpha k]} = \infty \right] & =
  \frac{\frac{1}{k}\sum_{i=[\alpha k]}^{k-1}\frac{q_i}{p_i} \cdots
    \frac{q_{[\alpha k] +1}}{p_{[\alpha k] +1}}}
       {\frac{1}{k}\sum_{i=[\alpha k]}^{\infty}\frac{q_i}{p_i} \cdots
         \frac{q_{[\alpha k]+1}}{p_{[\alpha k]+1}}}\\
& \underset{k \to \infty}{\longrightarrow} 
\frac{\int_{\alpha}^1 \left( \frac{\alpha}{t} \right)^8 dt}
     {\int_{\alpha}^{\infty} \left( \frac{\alpha}{t} \right)^8 dt} =
     1 - \alpha^7
\end{align*}
ce qui donne le r\'esultat voulu en prenant $\alpha$ assez petit.
\end{proof}

\begin{proposition}
\label{limlabspine}
Soit $Z$ un processus de Bessel de dimension $9$ issu de $0$. On a:
\[
 \left( \frac{1}{\sqrt{n}} Y_{[nt]} \right)_{t \geqslant 0} \underset{n \to
   \infty}{\longrightarrow} \left( Z_{\frac{2}{3} t} \right)_{t \geqslant 0}
\]
au sens de la convergence en loi dans l'espace $D(\mathbb{R}_+,\mathbb{R}_+)$.
\end{proposition}
\begin{proof}
La convergence voulue est une application directe d'un
r\'esultat plus g\'en\'eral de Lamperti \cite{Lam} que nous rappelons
maintenant. Soit $(X_n)_{n \geqslant 0}$ une cha\^ine de Markov homog\`ene sur
$\mathbb{R}_+$ v\'erifiant:
\begin{enumerate}
\item pour tout $K > 0$ on a de fa\c{c}on uniforme en
  $x \in \mathbb{R}_+$
\[\lim_{n \rightarrow \infty} \frac{1}{n} \sum_{i = 0}^{n-1}
\mathbb{P} \left( X_i \leqslant K \, \middle| \, X_0 = x \right) = 0;\]
\item pour tout $k \in \mathbb{N}$ les moments suivants existent et
  sont born\'es en tant que fonctions de la variable $x \in \mathbb{R}_+$
\[ m_k (x) = \mathbb{E} \left[ (X_{n+1} - X_n)^k | X_n = x \right];\]
\item il existe $\beta > 0$ et $\alpha > - \beta /2$ tels que
\begin{align*}
\lim_{x \rightarrow \infty} m_2(x) & = \beta,\\
\lim_{x \rightarrow \infty} x \, m_1(x) & = \alpha.
\end{align*}
\end{enumerate}
D\'efinissons le processus $(x_t^{(n)})_{t \in \mathbb{R}_+}$ par $x_t^{(n)} = n^{-1/2} X_i$
si $t = \frac{i}{n}$, $i = 0,1,2, \ldots$, ses valeurs en $t \in \mathbb{R}_+$
quelconque \'etant d\'etermin\'ees par interpolation lin\'eaire. Alors
le th\'eor\`eme de Lamperti \'enonce que
$(x_t^{(n)})_{t \in \mathbb{R}_+}$ converge en loi vers la diffusion
$(x_t)_{t \in \mathbb{R}_+}$ de g\'en\'erateur
\[ L = \frac{\alpha}{x} \frac{\mathrm{d}}{\mathrm{d}x} + \frac{\beta}{2}
\frac{\mathrm{d}^2}{\mathrm{d}x^2}.\]

Dans notre cas la cha\^ine de Markov \`a consid\'erer est $\widetilde{Y}$, de
matrice de transition $\widetilde{\Pi}$ pour laquelle $\widetilde{\Pi}
(x,y) = \Pi ([x],[y])$ si $y = x+1,x-1$ ou $x$.
La propri\'et\'e 1. découle alors facilement du lemme
\ref{labspine}. La propriété 2. est satisfaite trivialement.
De plus $p_n = \frac{1}{3} +
\frac{4}{3n} + {\it O} (n^{-2})$ et $q_n = \frac{1}{3} -
\frac{4}{3n} + {\it O} (n^{-2})$ (\cite{CD}, lemme 5.5) et donc:
\begin{align*}
\lim_{x \rightarrow \infty} m_2(x) & = \frac{2}{3}\\
\lim_{x \rightarrow \infty} x \, m_1(x) & = \frac{8}{3}.
\end{align*}
d'o\`u la propriété 3. avec $\alpha = 8/3$ et $\beta = 2/3$.

La cha\^ine $Y$ chang\'ee d'\'echelle converge donc en loi vers la
diffusion de g\'en\'erateur 
\[L = \frac{2}{3} \left(\frac{4}{x} \frac{\mathrm{d}}{\mathrm{d}x}
+ \frac{1}{2} \frac{\mathrm{d}^2}{\mathrm{d}x^2} \right),\]
d'o\`u le r\'esultat voulu.
\end{proof}

\subsection{Propri\'et\'es asymptotiques des petites \'etiquettes dans
les arbres}
\label{petitesetiquettes}

Le corollaire \ref{pbcle} permet de ramener la preuve du théorème \ref{egalite}
à la convergence des $\mu_n$-mesures de certaines boules dans l'espace des arbres
vers la $\mu$-mesure de ces mêmes boules, à ceci près qu'il faut contrôler l'erreur commise en écartant
les arbres qui sont dans $\Omega_S(R)$. Les resultats qui suivent ont pour but de contrôler cette erreur.
On fixe dans cette partie $R$ et $\varepsilon > 0$, et on note
$\Omega_S = \Omega_S(R)$.

\begin{lemme}
\label{seuil}
Il existe un entier $S^{\star}> 0$ tel que
$\marbres[] \left( \Omega_S \right) < \varepsilon$ pour tout $S > S^{\star}$.
\end{lemme}
\begin{proof}
Soit  $\Omega = \bigcap_{S=1}^{\infty} \Omega_S$. Si $\omega \in \Omega$,
$\omega$ a une infinit\'e de sommets dont les \'etiquettes appartiennent
\`a $\{1, \ldots, R+1 \}$, et donc il existe un entier $l \in\{1, \ldots,
R+1 \}$ tel que $N_l(\omega) = \infty$. Comme $\marbres[]$
est port\'ee par $\spinelab$ on a $\marbres[] ( \Omega)
= 0$ et donc $\marbres[] (\Omega_S) \rightarrow 0$ quand
$S \rightarrow \infty$.
\end{proof}

L'ingredient essentiel de la preuve du théorème \ref{egalite} est la proposition \ref{seuilN}
qui donne une estimation analogue à celle du lemme \ref{seuil} lorsque $\mu$ est remplacée par $\mu_n$,
\emph{avec uniformité} en $n$. Pour établir cette estimation, nous devrons d'abord majorer la probabilité
qu'à la génération $S$ existe un sommet d'étiquette plus petite que $S^{\alpha}$, o\`u $\alpha < 1/2$ est fixé:
c'est l'objet du lemme \ref{estimlabels} ci-dessous. Nous commençons par un lemme préliminaire facile:

\begin{lemme}
\label{taille}
Soit $S > 0$.
Il existe deux entiers positifs $N_1(S)$ et $K_{\varepsilon}(S)$ tels que pour tout $n > N_1(S)$:
\[ \marbres \left( \omega : \, \left|\boule[S]{\arbres}(\omega)\right|
  > K_{\varepsilon}(S)  \right) < \varepsilon .\]
\end{lemme}
\begin{proof}
C'est une cons\'equence directe de la convergence des mesures $\marbres$ vers
$\marbres[]$. En effet, comme
$|\boule[S]{\arbres}(\omega)| < \infty$ pour tout arbre $\omega$,
on peut choisir $K_{\varepsilon}(S)$ assez grand tel que
\[ \marbres[] \left( \omega : \, \left|\boule[S]{\arbres}(\omega) \right|
> K_{\varepsilon}(S)  \right) < \varepsilon .\]
La convergence de $\marbres$ vers $\marbres[]$ donne alors l'existence de
$N_1(S)$ tel que l'in\'egalit\'e du lemme soit valable pour tout $n > N_1(S)$.
\end{proof}

Il existe un nombre fini d'arbres bien \'etiquet\'es de hauteur exactement
$S$ et ayant au plus $K_{\varepsilon}(S)$ ar\^etes. On note
$M_{\varepsilon}(S)$ ce nombre.

\bigskip

Pour tout entier $S > 0$ et tout $\alpha \in \left[0, \frac{1}{2} \right[$ posons:
\[A_{\alpha}(S) = \left\{ \omega \in \arbres : \,
\omega \text{ a un sommet \`a la g\'en\'eration $S$
d'\'etiquette } \leqslant S^{\alpha}
\right\}.\]

\begin{lemme}
\label{estimlabels}
Soit $\alpha < \frac{1}{2}$. Pour tout entier $S$ assez grand, il existe
$N_2(S)$ tel que pour tout $n > N_2(S)$ on a:
\[\marbres \left( A_{\alpha}(S) \right)< \varepsilon.\]
\end{lemme}
\begin{proof}
Remarquons dans un premier temps qu'il suffit de d\'emontrer
l'in\'egalit\'e du lemme avec $\marbres$ remplac\'ee par $\marbres[]$.
En effet l'ensemble $A_{\alpha}(S)$ est un ferm\'e de $\arbres$. Si l'on sait
que $\marbres[] \left( A_{\alpha}(S) \right)< \varepsilon$, la
convergence des mesures $\marbres$ vers $\marbres[]$ assure alors que
$\limsup_N \marbres \left( A_{\alpha}(S) \right) \leqslant \marbres[]
\left( A_{\alpha}(S) \right) < \varepsilon$.

Rappelons les notations $\rho^{(l)}$ et $\widehat{\rho}^{(l)}$ introduites dans la partie \ref{sec:indirect}.
Pour $H > 0$ et $l>0$ on a:
\begin{equation*}
\widehat{\rho}^{(l)} \left( \vphantom{\rho^l} h(\omega) > H \right) =
\frac{1}{w_l} \sum_{\substack{\omega \in \arbre[{}]^{(l)}\\ h(\omega) > H}} 12^{-|\omega|}
\leqslant
\frac{1}{w_l} \sum_{\substack{\omega \in \mathbf{T}^{(l)}\\ h(\omega) > H}} 12^{-|\omega|}
= \frac{1}{w_l} \rho^{(l)} \left( \vphantom{\rho^l} h(\omega) > H \right).
\end{equation*}
Et donc
\[
\widehat{\rho}^{(l)} \left( \vphantom{\rho^l} h(\omega) > H \right)
\leqslant \frac{2}{w_l} \mathbb{P}_{GW(1/2)} \left[\vphantom{\rho^l} h(\omega) > H \right],
\]
o\`u $\mathbb{P}_{GW(1/2)}$ d\'esigne la loi de l'arbre de Galton-Watson de loi de
reproduction g\'eom\'etrique de param\`etre $1/2$. D'apr\`es le
th\'eor\`eme 1 (page 19) de \cite{AN}:
\[
\lim_{H \rightarrow \infty} H \, \mathbb{P}_{GW(1/2)} \left[
  \vphantom{\rho^k} h(\omega) > H \right] = 1.
\]
En observant que $\frac{2}{w_l} \leqslant \frac{3}{2}$ pour tout $l > 0$, on voit qu'il existe $H_1 > 0$ tel que si $H > H_1$:
\[
\widehat{\rho}^{(l)} \left(\vphantom{\rho^l}  h(\omega) > H \right) \leqslant \frac{2}{H}.
\]
Soit $\eta \in \left]0, \frac{1}{2}\right[$. Rappelons que $g_S(\omega)$ d\'esigne l'ensemble
des sommets de g\'en\'eration $S$ de $\omega$ et que pour $k$
entier, $L_k$ et $R_k$ sont les sous-arbres de $\omega$ attach\'es
respectivement \`a gauche et \`a droite du site de g\'en\'eration $k$
de la colonne de $\omega$. Pour tout $S > (1- \eta)^{-1} H_1$
on a d'apr\`es la deuxième partie du théorème \ref{descmu} et l'inégalité précédente:
\[
\marbres[] \left[ g_S(\omega) \cap \bigcup_{1 \leqslant k \leqslant
    [\eta S] - 1} (L_k \cup R_k) \neq \emptyset \right] \leqslant 2
\sum_{k=1} ^{[\eta S] - 1} \frac{2}{S-k} \leqslant 4 \frac{\eta}{1
  - \eta} \leqslant 8 \eta.
\]
Donc:
\[
\marbres[] \left( A_{\alpha}(S) \right) \leqslant 8 \eta +
\marbres[] \left( \exists s \in g_S(\omega) \cap \bigcup_{k = [\eta S]}^{S}
\left( L_k \cup R_k \right) \, : \, \ell(s) \leqslant S^{\alpha} \right).
\]
En appliquant la propriété de Markov de la chaîne $Y$ à l'instant $[\eta S]$,
puis en utilisant la proposition \ref{limlabspine} et le lemme \ref{labspine}, on trouve $\delta > 0$
et $S_1$ tels que pour tout $S > S_1$ on ait:
\[
\marbres[] \left( Y_k \geqslant [\delta \sqrt{S}] , \, \forall k \geqslant
  [\eta S] \right) \geqslant 1 - \eta.
\]
Il vient
\begin{equation}
\label{Aalphaintermediaire}
\marbres[] \left( A_{\alpha}(S) \right) \leqslant 9 \eta
+ \marbres[] \left( \left\{ \exists s \in g_S(\omega) \cap \bigcup_{k = [\eta S]}^{S}
\left( L_k \cup R_k \right) \, : \, \ell(s) \leqslant S^{\alpha} \right\} \cap \left\{
\forall k \geqslant [\eta S], \, Y_k \geqslant [\delta \sqrt{S}] \right\}
\right)
\end{equation}
Fixons-nous une famille $(y_k)_{[\eta S] \leqslant k \leqslant S}$ telle
que $y_k \geqslant [\delta \sqrt{S}]$ pour tout $k$. D'après le théorème \ref{descmu},
\begin{align}
\label{conditionel}
\marbres[] & \left( \exists s \in g_S(\omega) \cap \bigcup_{k = [\eta S]}^{S}
\left( L_k \cup R_k \right) \, : \, \ell(s) \leqslant S^{\alpha} \middle|
Y_k = y_k, \, [\eta S] \leqslant k \leqslant S \right) \notag \\
& \leqslant
2 \sum_{k = [\eta S]}^{S}
\widehat{\rho}^{(y_k)} \left( \vphantom{\sqrt{S}} \exists s \in g_{S-k} (\omega) :
  \, \ell(s) \leqslant S^{\alpha} \right)
 = 2 \sum_{k=0}^{S - [\eta S]}
\widehat{\rho}^{(y_{S-k})} \left(\vphantom{\sqrt{S}} \exists s \in g_{k} (\omega) :
  \, \ell(s) \leqslant S^{\alpha} \right).
\end{align}
Si $0 \leqslant k \leqslant S - [\eta S]$, on a:
\begin{equation*}
\begin{split}
\widehat{\rho}^{(y_{S-k})} \left(\vphantom{\sqrt{S}} \exists s \in g_{k} (\omega) :
  \, \ell(s) \leqslant S^{\alpha} \right) &  \leqslant
\widehat{\rho}^{(y_{S-k})} \left(\vphantom{\sqrt{S}}
\inf_{s \in \omega} \ell(s)  \leqslant S^{\alpha} \right)  \\
& = \frac{1}{w_{y_{S-k}}} \sum_{\substack{\omega \in \arbre[{}]^{(y_{S-k})}\\ \inf_{s \in \omega}
  \ell(s)  \leqslant S^{\alpha}}} 12^{- |\omega|}
 = \frac{1}{w_{y_{S-k}}} \sum_{\substack{\omega \in
  \mathbf{T}^{(y_{S-k})}\\ 0 < \inf_{s \in \omega}
  \ell(s)  \leqslant S^{\alpha} }} 12^{- |\omega|}\\
& = \frac{1}{w_{y_{S-k}}} \rho^{(y_{S-k})} \left(\vphantom{\sqrt{S}}
 0 < \inf_{s \in \omega} \ell(s) \leqslant S^{\alpha} \right) .\\
\end{split}
\end{equation*}
Or
\[
\rho^{(y_{S-k})} \left(\vphantom{\sqrt{S}}
\inf_{s \in \omega} \ell(s) > 0 \right) = w_{y_{S-k}} 
\]
et
\begin{align*}
\rho^{(y_{S-k})} \left(\vphantom{\sqrt{S}}
\inf_{s \in \omega} \ell(s) > S^{\alpha} \right) 
\geqslant \rho^{(y_{S-k})} \left(\vphantom{\sqrt{S}}
\inf_{s \in \omega} \ell(s) > [S^{\alpha}] \right) 
& = \rho^{(y_{S-k} - [S^{\alpha}])} \left(\vphantom{\sqrt{S}} \inf_{s
    \in \omega} \ell(s) > 0 \right) \notag \\
& = w_{y_{S-k} - [S^{\alpha}]}. 
\end{align*}
On a donc 
\[
\widehat{\rho}^{(y_{S-k})} \left(\vphantom{\sqrt{S}} \exists s \in g_{k} (\omega) :
  \, \ell(s) \leqslant S^{\alpha} \right)  \leqslant
\frac{1}{w_{y_{S-k}}} \left( w_{y_{S-k}} - w_{y_{S-k} -  [S^{\alpha}]}
\right) = 1 - \frac{w_{y_{S-k} - [S^{\alpha}]}}{w_{y_{S-k}}}.
\]
Comme $y_{S-k} \geqslant [ \delta \sqrt{S}]$, un d\'eveloppement limit\'e nous donne
\[
1 - \frac{w_{y_{S-k} - [S^{\alpha}]}}{w_{y_{S-k}}} =
4 S^{\alpha} y_{S-k}^{-3} + o\left( S^{\alpha} y_{S-k}^{-3} \right)
 \leqslant \frac{4}{\delta^3} S^{\alpha - \frac{3}{2}} +  o\left( S^{\alpha - \frac{3}{2}}  \right),\]
et donc le terme de gauche de \eqref{conditionel} est majoré par
$\frac{8}{\delta^3} S^{\alpha -1/2} + o \left(S^{\alpha - 1/2} \right)$.
En revenant à \eqref{Aalphaintermediaire}, on a
\[
\marbres[] \left( A_{\alpha}(S) \right) \leqslant 9 \eta +
\frac{9}{\delta^3} S^{\alpha - \frac{1}{2}} +  o\left( S^{\alpha - \frac{1}{2}} \right).
\]
Finalement, puisque $\alpha < \frac{1}{2}$, on trouve que $\mu
(A_{\alpha}(S)) < 10 \eta$ d\`es que $S$ est assez grand, ce qui termine la preuve.
\end{proof}

Ce dernier lemme va nous permettre de d\'emontrer la version uniforme en $n$ du lemme \ref{seuil}:
\begin{proposition}
\label{seuilN}
Pour tout entier $S$ assez grand, il existe $N(S)$ tel que pour tout
$n > N(S)$ on a:
\[\marbres(\Omega_S)< \varepsilon .\]
\end{proposition}
\begin{proof}
Dans toute la preuve on se fixe $\alpha \in \left]\frac{1}{3},
\frac{1}{2} \right[$. Si $S > 0$, le lemme \ref{taille} donne
$K_{\varepsilon}(S) > 0$ et $N_1(S)> 0$ tels que si $n > N_1(S)$ on a 
$\marbres \left( \omega : \, |\boule[S]{\arbres}(\omega)| > K_{\varepsilon}(S)  \right) <
\varepsilon $. Rappelons aussi que le nombre d'arbres bien \'etiquet\'es
de hauteur $S$ et de taille inf\'erieure \`a $K_{\varepsilon}(S)$ est
$M_{\varepsilon}(S)$. Le lemme \ref{estimlabels} montre aussi que pour
$S$ assez grand, il existe $N_2(S)$ tel que pour tout $n > N_2(S)$ on a
$\marbres(A_{\alpha}(S)) < \varepsilon$. On a donc pour tout $S$ assez grand et pour tout $n >
N_1(S) \vee N_2(S)$:
\begin{align}
\marbres (\Omega_S)  & = \sum_{\substack{\omega^{\star} \notin A_{\alpha}(S)\\
|\omega^{\star}| \leqslant K_{\varepsilon}(S), \, h(\omega^{\star}) = S}} \mu_n \left( \{\omega 
  : \, \boule[S]{\arbres}(\omega) = \omega^{\star} \} \cap
\Omega_S \right) \notag\\
& \qquad + \marbres(A_{\alpha}(S)) +  \marbres
\left( \omega: \, |\boule[S]{\arbres}(\omega)| > K_{\varepsilon}(S)
\right) \notag \\
& \leqslant 2 \varepsilon + \sum_{\substack{\omega^{\star} \notin A_{\alpha}(S)\\
|\omega^{\star}| \leqslant K_{\varepsilon}(S), \, h(\omega^{\star}) = S}} \mu_n \left( \{\omega 
  : \, \boule[S]{\arbres}(\omega) = \omega^{\star} \} \cap
\Omega_S  \right).
\label{munosr}
\end{align}

Fixons-nous un arbre $\omega^{\star} \notin A_{\alpha}(S)$, de taille
inf\'erieure \`a $K_{\varepsilon}(S)$, de hauteur $S$ et
dont les sommets de g\'en\'eration
$S$, au nombre de $k$, sont d'\'etiquettes $l_1, \ldots, l_k$. En
d\'ecomposant suivant les sous-arbres de $\omega$ ayant pour racines les
sommets de la g\'en\'eration $S$, de fa\c{c}on analogue \`a
\eqref{boulen}, on obtient l'in\'egalit\'e:
\begin{equation}
\label{omegastar}
\marbres \left( \{\omega : \, \boule[S]{\arbres}(\omega) = \omega^{\star} \} \cap
\Omega_S \right) \leqslant
\frac{1}{D_n} \, \sum_{n_1 + \cdots + n_k = n - |\omega^{\star}|} \, \sum_{i = 1}^k
D_{n_i}^{(l_i)}(R) \prod_{j \neq i} D_{n_j}^{(l_j)}
\end{equation}
o\`u $D_n^{(l)}(R)$ est le nombre d'arbres de $\arbre[n]^{(l)}$
ayant au moins un sommet d'\'etiquette $R$.
Comme $\omega \notin A_{\alpha}(S)$, si $S$ est assez grand, on a
$l_i > S ^{\alpha} > R$, et donc $D_{n_i}^{(l_i)}(R) = D_{n_i}^{(l_i)} - D_{n_i}^{(l_i - R)}$
pour tout $i = 1, \ldots , k$, et (\ref{omegastar}) devient:
\begin{equation*}
\begin{split}
\mu_n  & \left( \{\omega : \, \boule[S]{\arbres}(\omega) = \omega^{\star} \} \cap
\Omega_S \right) \\
& \leqslant
\frac{1}{D_n} \, \sum_{n_1 + \cdots + n_k = n - |\omega^{\star}|} \, \sum_{i = 1}^k
(D_{n_i}^{(l_i)} - D_{n_i}^{(l_i - R)}) \prod_{j \neq i} D_{n_j}^{(l_j)}\\
& =
\frac{k}{D_n} \, \sum_{n_1 + \cdots + n_k = n - |\omega^{\star}|} \prod_{j=1}^k
D_{n_j}^{(l_j)}
- \frac{1}{D_n} \,  \sum_{i = 1}^k \, \sum_{n_1 + \cdots + n_k = n -
  |\omega^{\star}|} \, D_{n_i}^{(l_i - R)} \prod_{j \neq i} D_{n_j}^{(l_j)}.
\end{split}
\end{equation*}

La convergence $\marbres \left(\omega : \boule[S]{\arbres} (\omega) =
\omega^{\star} \right) \rightarrow \marbres[] \left(\omega : \boule[S]{\arbres} (\omega) =
\omega^{\star} \right)$  quand $n \rightarrow \infty$ (th\'eor\`eme
\ref{cvarbres}) ainsi que les identit\'es (\ref{boulen}) et
(\ref{boule}) permettent de trouver $N(\omega^{\star},S)$ tel que pour tout $n
> N(\omega^{\star},S)$:
\begin{equation*}
\frac{1}{D_n} \,  \sum_{n_1 + \cdots + n_k = n - |\omega^{\star}|} \prod_{j=1}^k
D_{n_j}^{(l_j)}
\leqslant
12^{- |\omega^{\star}|} \, \sum_{t = 1}^k d_{l_t} \prod_{s \neq t} w_{l_s} \, +
\, \frac{\varepsilon}{K_{\varepsilon}(S)M_{\varepsilon}(S)}
\end{equation*}
et, pour tout $i = 1 ,\ldots, k$:
\begin{equation*}
\begin{split}
 \frac{1}{D_n}\, & \sum_{n_1 + \cdots + n_k = n - |\omega^{\star}|} \,
D_{n_i}^{(l_i - R)} \prod_{j \neq i} D_{n_j}^{(l_j)} \\
& \geqslant
12^{- |\omega^{\star}|} \Big( d_{l_i - R} \prod_{j \neq i} w_{l_j} + \sum_{t \neq
  i} d_{l_t} w_{l_i - R} \prod_{j \neq t, i} w_{l_j} \Big) -
\frac{\varepsilon}{K_{\varepsilon}(S)M_{\varepsilon}(S)}.
\end{split}
\end{equation*}
On a donc pour tout $n > N(\omega^{\star},S)$:
\begin{align}
\marbres & \left( \{\omega : \, \boule[S]{\arbres}(\omega) = \omega^{\star} \} \cap
\Omega_S \right) \notag \\
%%%%%%%%%%
 & \leqslant \frac{2\varepsilon}{M_{\varepsilon}(S)}
+ k 12^{- |\omega^{\star}|} \, \sum_{t = 1}^k d_{l_t} \prod_{s \neq t} w_{l_s}
%%%%%%%%%%
- 12^{- |\omega^{\star}|} \sum_{i=1}^k \Big( d_{l_i - R} \prod_{j \neq i} w_{l_j} + \sum_{t \neq
  i} d_{l_t} w_{l_i - R} \prod_{j \neq t, i} w_{l_j} \Big) \notag \\
%%%%%%%%%%%%%%%%%%%%%%%%%
 & \leqslant \frac{2\varepsilon}{M_{\varepsilon}(S)}
+ 12^{- |\omega^{\star}|} \, \sum_{t = 1}^k (d_{l_t} - d_{l_t -R}) \prod_{s \neq
  t} w_{l_s}
+ 12^{- |\omega^{\star}|} \, \sum_{t = 1}^k d_{l_t}\left( \sum_{i \neq
    t} (w_{l_i} - w_{l_i - R}) \prod_{s \neq t,i} w_{l_s} \right).\label{muncapomega}
\end{align}
Notons:
\begin{align*}
d(\omega^{\star}) & = \max_{i= 1 \ldots k} \left( 1 - \frac{d_{l_i - R}}{d_{l_i}}
\right),\\
w(\omega^{\star}) & =  \max_{i= 1 \ldots k} \left( 1 - \frac{w_{l_i - R}}{w_{l_i}}
\right).
\end{align*}
L'in\'egalit\'e \eqref{muncapomega} devient, pour tout $n > N(\omega^{\star},S)$:
\begin{align}
\marbres & \left( \{\omega : \, \boule[S]{\arbres}(\omega) = \omega^{\star} \} \cap
\Omega_S \right) \notag \\
%%%%%%%%%%
 & \leqslant \frac{2\varepsilon}{M_{\varepsilon}(S)}
+ d(\omega^{\star}) \, 12^{- |\omega^{\star}|} \, \sum_{t = 1}^k d_{l_t} \prod_{s \neq
  t} w_{l_s}
%%%%%%%%%%
+ k w(\omega^{\star}) \, 12^{- |\omega^{\star}|} \, \sum_{t = 1}^k d_{l_t} \prod_{s \neq
  t} w_{l_s} \notag \\
& = \frac{2\varepsilon}{M_{\varepsilon}(S)}
+ \left( d(\omega^{\star}) + k . w(\omega^{\star}) \right) \marbres[] \left( \omega :
\boule[s]{\arbres} (\omega) = \omega^{\star} \right)
\label{omegastarfinal}
\end{align}
en utilisant \eqref{boule} pour la dernière inégalité.

\bigskip

Posons maintenant
$N^{\star}(S) = \max_{|\omega^{\star}| \leqslant K_{\varepsilon}(S)} N(\omega^{\star},S) \vee
N_1(S) \vee N(S)$. Si $S$ est assez grand, pour $n > N^{\star}(S)$ on a en
utilisant (\ref{munosr}) et (\ref{omegastarfinal}):
\begin{equation*}
\marbres(\Omega_S) \leqslant 4 \varepsilon
+ \sum_{\substack{\omega^{\star} \notin A_{\alpha}(S)\\ |\omega^{\star}| \leqslant
  K_{\varepsilon}(S), \, h(\omega^{\star}) = S}}
\left( d(\omega^{\star}) + |g_S(\omega^{\star})| . w(\omega^{\star}) \right) \marbres[] \left(\omega :
\boule[S]{\arbres} (\omega) = \omega^{\star} \right).
\end{equation*}

Un d\'eveloppement limit\'e donne, pour $\omega^{\star} \notin A_{\alpha}(S)$, $w(\omega^{\star}) \leqslant
4 (5R - 3) S^{-3 \alpha} + o\left(S^{-3 \alpha} \right)$.
De plus, $\sup_{\omega^{\star} \notin A_{\alpha}(S)} d(\omega^{\star})
\rightarrow 0$ quand $S \rightarrow \infty$.
Tout cela permet de trouver $S^{\star}$ tel que si $S > S^{\star}$ et $n> N^{\star}(S)$:
\begin{align}
\marbres(\Omega_S) & \leqslant 4 \varepsilon +
\sum_{\substack{\omega^{\star} \notin A_{\alpha}(S)\\  |\omega^{\star}| \leqslant K(S), \, h(\omega^{\star}) = S}}
\left( \varepsilon + |g_S(\omega^{\star})| . 20R S^{-3 \alpha} \right) \marbres[] \left(\omega :
\boule[S]{\arbres}(\omega) = \omega^{\star} \right) \notag \\
& \leqslant 5 \varepsilon +
20R S^{-3 \alpha} \sum_{\substack{\omega^{\star} \notin A_{\alpha}(S)\\ |\omega^{\star}| \leqslant
  K(S), \, h(\omega^{\star})=S}} |g_S(\omega^{\star})|  \, \marbres[] \left(\omega :
\boule[S]{\arbres}(\omega) = \omega^{\star} \right) \notag \\
& \leqslant 5 \varepsilon +
20R S^{-3 \alpha} \mathbb{E}_{\marbres[]}
\left[|g_S(\omega)| \right]. \label{munomegas}
\end{align}

La description de $\marbres[]$ donn\'ee dans
le th\'eor\`eme \ref{descmu} permet d'estimer
$\mathbb{E}_{\marbres[]} \left[|g_S(\omega)| \right]$. En effet, si $H > 0$, on a pour tout entier $k \geqslant 1$
\[
\mathbb{E}_{\widehat{\rho}^{(k)}} \left[ |g_H(\omega)| \right] \leqslant
\frac{1}{w_k} \mathbb{E}_{\rho^{(k)}} \left[ |g_H(\omega)| \right] =
\frac{2}{w_k} \mathbb{E}_{GW(1/2)} \left[ |g_H(\omega)| \right] =
\frac{2}{w_k} \leqslant 2.
\]
Donc:
\[
\mathbb{E}_{\marbres[]} \left[ |g_S(\omega)| \right] \leqslant 4S +1.
\]
En utilisant \eqref{munomegas}, et en rappelant que $\alpha > \frac{1}{3}$, on voit que pour tout $S$ assez grand
et pour $n > N^{\star}(S)$ on a
\[
\mu_n(\omega_S) \leqslant 6 \varepsilon.
\]
Quitte à remplacer $\varepsilon$ par $\varepsilon/6$ dans la preuve ci-dessus, cela donne la majoration annoncée.
\end{proof}

\subsection{Preuve du r\'esultat principal}
\label{deroulement}

Dans toute la preuve, $q^{\star} \in \kcartes$ et $R > 0$
sont fix\'es et $\Omega_S = \Omega_S(R)$ comme ci-dessus. Le théorème \ref{egalite}
est \'equivalent \`a la convergence suivante:
\begin{equation}
\label{conv}
\marbres \Big( \omega :
\boule{\kcartes}(\bij(\omega)) = \boule{\kcartes}(q^{\star}) \Big)
\underset{n \to \infty}{\longrightarrow}
\marbres[] \Big( \omega :
\boule{\kcartes}(\bij(\omega)) = \boule{\kcartes}(q^{\star}) \Big).
\end{equation}

Il faut dans un premier temps reformuler le probl\`eme en termes
d'arbres. Comme $q^{\star} \in \kcartes$, il existe une
quadrangulation finie  $q_0 \in \carte[]$ telle que $\dist{\kcartes}(q_0,q^{\star}) <
\frac{1}{R+1}$, donc $\boule{\kcartes} (q_0) = \boule{\kcartes} (q^{\star})$.
De plus $q_0$ \'etant une quadrangulation
finie, il existe $\omega_0 \in \arbre[]$ tel que $\bij(\omega_0) = q_0$,
et donc $\boule{\kcartes} (\bij(\omega_0)) = \boule{\kcartes} (q^{\star})$.
(\ref{conv}) devient alors :
\begin{equation}
\label{convbis}
\marbres \Big( \omega :
\boule{\kcartes} (\bij(\omega)) = \boule{\kcartes}(\bij(\omega_0)) \Big)
\underset{n \to \infty}{\longrightarrow} \marbres[] \Big( \omega :
\boule{\kcartes} (\bij(\omega)) = \boule{\kcartes} (\bij(\omega_0)) \Big).
\end{equation}
On fixe $\varepsilon > 0$ dans toute la suite.

\bigskip

Le probl\`eme est maintenant de caract\'eriser les arbres
$\omega$ pour lesquels la quadrangulation associ\'ee a une boule de rayon
$R$ qui co\"incide avec celle de $\bij(\omega_0)$. La difficult\'e
vient du fait que deux arbres
proches dans $\arbres$ peuvent donner des quadrangulations
tr\`es diff\'erentes s'ils ont des sommets de petite \'etiquette dans des
g\'en\'erations \'elev\'ees. Ce probl\`eme se r\'esoud \`a l'aide de la proposition \ref{seuilN}.

Observons que $\omega_0$ est un arbre fini. Si $S_0$ d\'esigne sa hauteur,
$\omega_0$ n'a \'evidemment pas d'\'etiquette inf\'erieure ou \'egale \`a
$R+1$ strictement au dessus de la g\'en\'eration $S_0$. Grâce au lemme \ref{seuil} et à la proposition
\ref{seuilN}, on peut choisir $S_1 > S_0$ tel que pour tout $S \geqslant S_1$ on a $\mu(\Omega_S) < \varepsilon$
et $\mu_n(\Omega_S) < \varepsilon$ pour tout $n \geqslant N(S)$.

Soit $S > S_1$ et $(\omega_i)_{i \in I}$ la famille du corollaire \ref{pbcle},
telle que pour tout $\omega \in \spinelab \cap \Omega_S^c$ l'égalité $\boule{\kcartes} (\bij(\omega)) =
\boule{\kcartes} (\bij(\omega_0))$ a lieu si et seulement si il existe $i \in I$ tel que
$\boule[S]{\arbres} (\omega) = \boule[S]{\arbres} (\omega_i)$. En notant $A
\vartriangle B$ la diff\'erence sym\'etrique des ensembles $A$ et $B$ on a
\begin{equation*}
\begin{split}
\marbres[] & \left( \left\{ \omega \in \arbres : \,
\boule{\kcartes} (\bij(\omega)) = \boule{\kcartes} (\bij(\omega_0)) \right\}
\vartriangle
\bigcup_{i\in I} \left\{ \omega \in \arbres : \, 
\boule[S]{\arbres} (\omega) = \boule[S]{\arbres} (\omega_i) \right\}
\right)\\
& \leqslant \marbres[] \left( \Omega_S \right) < \varepsilon.
\end{split}
\end{equation*}
Il en d\'ecoule que:
\begin{equation*}
\begin{split}
\Bigg| & \marbres \Big( \omega : \,
\boule{\kcartes} (\bij(\omega)) = \boule{\kcartes} (\bij(\omega_0)) \Big)
%%%%
- \marbres[] \Big( \omega : \,
\boule{\kcartes} (\bij(\omega)) = \boule{\kcartes} (\bij(\omega_0)) 
\Big) \Bigg|\\
%%%%%%%%%%%%%%%%%%%%%%%%%%%%%%%%%%%%%%%%%%%%%%%%%%%%%%%%%%%%%%%%%%%%%%%
& \leqslant  \Bigg| \marbres[] \Big(
\bigcup_{i\in I} \left\{ \omega : \,
\boule[S]{\arbres} (\omega) = \boule[S]{\arbres} (\omega_i) \right\} \Big)
%%%%
- \marbres \Big(\omega : \,
\boule{\kcartes} (\bij(\omega)) = \boule{\kcartes}
(\bij(\omega_0)) \Big) \Bigg|\\
%%%%%%%%%%%%%%%%%%%%%%%%%%%%%%%%%
 &  \qquad + \Bigg| \marbres[] \Big(
\bigcup_{i\in I} \left\{ \omega : \,
\boule[S]{\arbres} (\omega) = \boule[S]{\arbres} (\omega_i) \right\} \Big)
%%%%
- \marbres[] \Big( \omega : \,
\boule{\kcartes} ( \bij(\omega)) = \boule{\kcartes}
(\bij(\omega_0)) \Big) \Bigg|\\
%%%%%%%%%%%%%%%%%%%%%%%%%%%%%%%%%%%%%%%%%%%%%%%%%%%%%%%%%%%%%%%%%%%%%%%%
 & \leqslant \Bigg| \marbres[] \Big(
\bigcup_{i\in I} \left\{ \omega : \,
\boule[S]{\arbres} (\omega) = \boule[S]{\arbres} (\omega_i) \right\} \Big)
%%%%
- \marbres \Big( \omega : \,
\boule{\kcartes} (\bij(\omega)) = \boule{\kcartes} (\bij(\omega_0)) 
\Big) \Bigg| + \varepsilon .
\end{split}
\end{equation*}
De plus $\bigcup_{i\in I} \left\{ \omega \in \arbres : \, 
\boule[S]{\arbres} (\omega) = \boule[S]{\arbres} (\omega_i) \right\}$
est ouvert et ferm\'e dans $\overline{\mathbb{T}}$. Donc,
\begin{equation*}
\marbres \left( \bigcup_{i\in I} \left\{
\omega : \, \boule[S]{\arbres} (\omega) 
= \boule[S]{\arbres}(\omega_i) \right\} \right)
\underset{n \to \infty}{\longrightarrow}
\marbres[] \left( \bigcup_{i\in I} \left\{
\omega : \, \boule[S]{\arbres} (\omega) 
= \boule[S]{\arbres} (\omega_i) \right\} \right).
\end{equation*}
On obtient ainsi que pour $S > S_1$ il existe $N(S) > 0$ tel que pour tout $n > N(S)$:
\begin{equation*}
\begin{split}
\Bigg| & \marbres \Big( \omega : \,
\boule{\kcartes} (\bij(\omega)) = \boule{\kcartes} (\bij(\omega_0)) \Big)
%%%%
- \marbres[] \Big( \omega : \,
\boule{\kcartes} (\bij(\omega)) = \boule{\kcartes} (\bij(\omega_0)) \Big) \Bigg|\\
%%%%%%%%%%%%%%%%%%%%%%%%%%%%%%%%%%%%%%%%%%%%%%%%%%%%%%%%%%%%%%%%%%%%%%%
 & \leqslant \Bigg| \marbres \Big( \bigcup_{i\in I} \left\{ \omega : \,
\boule[S]{\arbres} (\omega) = \boule[S]{\arbres} (\omega_i) \right\} \Big)
%%%%
- \marbres \Big( \omega : \,
\boule{\kcartes} (\bij(\omega)) = \boule{\kcartes} (\bij(\omega_0)) \Big) \Bigg|
+ 2 \varepsilon \\
%%%%%%%%%%%%%%%%%%%%%%%%%%%%%%%%%%%%%%%%%%%%%%%%%%%%%%%%%%%%%%%%%%%%%%%
 & = \Bigg| \marbres \Big( \bigcup_{i\in I} \left\{ \omega : \,
\boule[S]{\arbres} (\omega) = \boule[S]{\arbres} (\omega_i) \right\} \cap
\Omega_S \Big)
%%%%
- \marbres \Big( \left\{ \omega : \,
\boule{\kcartes} (\bij(\omega)) = \boule{\kcartes} (\bij(\omega_0)) 
\right\} \cap \Omega_S \Big) \Bigg| \\
& \qquad + 2 \varepsilon .
\end{split}
\end{equation*}
d'apr\`es le choix de la famille $(\omega_i)_{i \in I}$.

\bigskip

Par ailleurs on a aussi $\mu_n(\Omega_S) < \varepsilon$ pour $n > N(S)$ et donc pour tout $n > N(S) \vee N_1(S)$,
\begin{equation*}
\Bigg| \marbres
\Big( \omega : \,
\boule{\kcartes} (\bij(\omega)) = \boule{\kcartes} (\bij(\omega_0)) \Big)
- \marbres[]
\Big( \omega : \,
\boule{\kcartes} (\bij(\omega)) = \boule{\kcartes} (\bij(\omega_0)) \Big) \Bigg|
\leqslant 3 \varepsilon
\end{equation*}
ce qui termine la preuve du théorème \ref{egalite}. \qed

\bigskip

\noindent {\bf Remerciements.} L'auteur souhaite remercier
Jean-Fran\c{c}ois Le Gall pour de nombreuses conversations sur le
pr\'esent travail.

\addcontentsline{toc}{section}{R\'ef\'erences}
\bibliographystyle{abbrv}
\def\cprime{$'$}

\end{document}